\documentclass[11pt]{amsart}
\usepackage{amsmath,amsthm,amsfonts,amssymb,enumerate}
\usepackage[normalem]{ulem}
\newtheorem{Theorem}{Theorem}[section]

\newtheorem{Lemma}[Theorem]{Lemma}
\newtheorem{Corollary}[Theorem]{Corollary}

\newtheorem{Remark}[Theorem]{Remark}

\newcommand{\bpf}[1][Proof]{{\noindent {\sc #1: }}}
\newcommand{\epf}{{{\hfill $\Box$ \smallskip}}}

\newcommand{\Domain}{O}
\newcommand{\ds}{\displaystyle}
\newcommand{\eps}{\varepsilon}

\newcommand{\E}{\mathsf{E}}
\newcommand{\Fc}{\mathcal{F}}
\newcommand{\Bc}{\mathcal{B}}
\newcommand{\Bb}{\mathbb{B}}

\newcommand{\Pp}{\mathsf{P}}

\newcommand{\R}{\mathbb{R}}
\newcommand{\N}{\mathbb{N}}

\newcommand{\ONE}{{\bf 1}}

\newcommand{\Tr}{\mathop{\rm Tr}}

\newcommand{\dist}{\mathop{\rm dist}}
\newcommand{\supp}{\mathop{\rm supp}}

\def\Swiech
{{\accent"13S}wie{\hbox{\kern -0.21em\lower
0.79ex\hbox{$\textfont1=\scriptfont1
\lhook$}}}ch}

\def\SWIECH
{{\accent"13S}WIE{\hbox{\kern -0.21em\lower
0.79ex\hbox{$\textfont1=\scriptfont1
\lhook$}}}CH}

\begin{document}
\baselineskip=17pt

\title[Scaling limits for conditional exit problems]{Scaling limits for conditional diffusion exit problems, Doob's $h$-transform, and asymptotics for nonlinear
elliptic equations. }

\date{}
\author{Yuri Bakhtin and Andrzej \Swiech}
\address{School of Mathematics\\ Georgia Institute of Technology\\
686 Cherry Street \\ Atlanta, GA, 30332-0160}
\thanks{Yuri Bakhtin was partially supported by the NSF CAREER Award DMS-0742424. Andrzej
\Swiech\ was partially supported by the NSF grant DMS-0856485.}


\begin{abstract} The goal of this paper is to supplement
the large deviation principle of the Freidlin--Wentzell theory on exit problems for diffusion processes 
with results of classical central limit theorem kind. We describe a class of situations where
conditioning on exit through unlikely locations leads to a Gaussian scaling limit for the exit distribution. Our results are based on Doob's $h$-transform
and new asymptotic convergence gradient estimates for elliptic nonlinear equations that allow to reduce the problem to the Levinson case.
 We devote a separate section to a rigorous and general discussion of $h$-transform.
 \end{abstract}



\maketitle

\section{Introduction}

The studies of dynamical systems under small white noise perturbations
have a long history. Although the limiting behavior of the resulting stochastic dynamics as the noise intensity vanishes
depends very much on the character of the system, in many cases the key features
can be described in terms of the celebrated Freidlin--Wentzell theory (FW), see~\cite{FW:MR722136}. 

Exit problems associated with these diffusions provide important information 
that can be used to study transitions of the system between various regions in the phase space. FW theory
provides the asymptotic description of the exit distributions at the level of large deviations estimates.
In several important situations besides finding the points where the exit distribution concentrates
as the noise intensity vanishes, it also produces exact exponential concentration rates.
A key notion in FW is that of quasi-potential, a rate of unlikelihood or cost of diffusing from one point to another.
One can often say that the exit distribution concentrates near minima of quasi-potential in the small noise limit, although
concrete results require more careful statements.

The goal of this paper is to supplement
the large deviation principle results of FW with a result of classical central limit theorem type, i.e., a Gaussian scaling limit.

Although the idea of obtaining such result is not new --- a possibility of such result in the large deviation framework
was hinted at in the original book~\cite{FW:MR722136} ---
no rigorous results are known to us. However, such results may have important consequences. For example, 
in~\cite{B:MR2731621,B:MR2800902,B:MR2802310}, the
scaling limits for exit distribution played the crucial role in the analysis of small noise limit for noisy
heteroclinic networks. In that context, in 
the exit problem near a saddle point of the drift vector field,
 the quasi-potential approach does not produce enough detail. At the same time the corresponding scaling limits capture 
the difference between asymptotically symmetric and asymmetric exit distributions responsible for Markov or non-Markov
behavior of the limiting process jumping between saddles along heteroclinic connections.

FW theory mainly concerns the situation where for any finite
time horizon the exit from a domain~$O$ gets extremely unlikely as the noise intensity~$\eps$ tends to zero.
Such situations arise if the unperturbed dynamical system has an attractor or several attractors in~$\Domain$.
Then the system spends
a long time in a small neighborhood of an attractor and rarely makes excursions away from it as the drift brings it back
since the noise is small. After a long time and multiple escape attempts, a larger fluctuation inevitably occurs
and makes one of these attempts succeed, so the system either exits the domain or approaches one
of the other attractors and the same scenario repeats there. The arising phenomenon of noise induced rare jumps between
attractors is often called metastability. 

In this setting, the reactive trajectories or transition paths usually travel
through a small neighborhood of a saddle critical point of the drift, and the location of the exit from that neighborhood may critically
impact the remainder of the transition path by influencing the choice of the next attractor visited by the
system. So, obtaining distributional scaling limits for such situations would be highly desirable as it may lead to
non-Markovian effects of asymmetrical decisions analogous to~\cite{B:MR2731621,B:MR2800902,B:MR2802310}.

We are not ready to make any mathematical claims concerning this important case. Instead we 
consider a simpler situation.
Although transition or exit without immediate return to the neighborhood of the attractor is a rare event, 
one can study the system conditioned on this rare event. The goal of this paper is to describe a class of situations where
conditioning on exit through unlikely locations leads to a Gaussian scaling limit for the exit distribution.

Let us briefly describe our program. 
 
The first step is to use Doob's $h$-transform that allows to claim that a diffusion conditioned on an exit event
is again a diffusion process with the original diffusion coefficient and a new drift that is obtained
from the original drift by a nonlinear correction that depends on~$\eps$. The $h$-transform is a well-known tool, 
see, e.g., a recent work~\cite{Nolen-Lu:2013arXiv1303.1744L} that uses $h$-transform to support the theory of transition
paths~\cite{EVE:MR2252154},\cite{MSVE06},\cite{EVE2010}. A general, rigorous, and detailed discussion of
$h$-transform is given in Section~\ref{sec:h-transform} playing the role of an appendix.

The second step is to establish
sufficiently fast convergence of the corrected drift
as $\eps\to 0$ to a limiting vector field.  The idea is that the convergence in question can often 
be
studied with the help of asymptotic analysis of stationary Hamilton--Jacobi--Bellman (HJB) equation with
positive viscosity~$\eps$. 

In fact, a large part of this paper is a discussion of a class of such situations where solutions 
of viscous HJB equations and their gradients converge in domains of regularity, as $\eps\to0$, to
viscosity solutions of the inviscid HJB equations
and, moreover, an expansion for the solution with the leading correction term of 
the order of~$\eps^2$ holds true in $C^1$ norm. Such expansions are known in $C^1$ for evolutionary
HJB equations, see~\cite{F1:MR0304045}, and in~$C^0$ for stationary ones with zero boundary condition, see~\cite{FS:MR833404}.  
No results on $C^1$ expansions in the stationary case with our boundary condition are known to us, and our result
seems to be the first one in this direction. We have to impose certain restrictions on the problem, namely,
we assume that the noise is additive and the boundary is flat. Our analysis 
is based
on PDE methods and involves elements of stochastic control. It is worth mentioning that it is the infinite time
horizon in the corresponding stochastic control problems that makes the analysis of stationary HJB equations
harder than that of evolutionary ones, where the time horizon is finite. Let us also mention the
connection between the FW quasi-potential, the rate function in the large deviation principle 
for the conditional exit point, and the viscosity solution of the stationary HJB equation (see~\cite{F2:MR512217}
for a more PDE and stochastic control viewpoint than~\cite{FW:MR722136}).

The third step in our program is to check that the limiting drift satisfies the Levinson conditions.
If it does, then the CLT for conditioned diffusion follows from the exit CLT for diffusions in 
the Levinson case
obtained in~\cite{B:MR2739004}.

We will describe the setting, explain our program in more detail, and formulate the main results in 
Section~\ref{sec:main-results}. In Section~\ref{sec:example} we show a very
simple example where our program can be carried out by explicit computations.
In Section~\ref{sec:PDE-main-results} we explain the connection with HJB equations
and state relevant results on asymptotics of solutions of HJB equations guaranteeing the CLT of Section~\ref{sec:main-results}.
Proofs of these results are given in Section~\ref{sec:PDE}. The auxiliary
Section~\ref{sec:h-transform} is devoted to Doob's $h$-transform.

\section{Main results}\label{sec:main-results}
Let us now be more precise. We begin with the deterministic dynamics in $\R^n$, $n\in\N$, defined by 
a~$C^2$~vector field $b(x)=(b^1(x),\ldots, b^n(x)), x\in\R^n$:
\[
\dot X(t)=b(X(t)).
\]
For $\eps>0$, we consider an elliptic stochastic perturbation of this system given by the following It\^o equation:
\begin{equation}
d X_\eps(t)=b(X_\eps(t))dt+\eps \sigma(X_\eps(t)) dW(t).
\label{eq:sde-perturbation}
\end{equation}
Here $W=(W^1,\ldots,W^m)$, $m\in\N$, is a standard $m$-dimensional Wiener process  defined
on a filtered probability space $(\Omega,\Fc, (\Fc_t)_{t\ge 0},\Pp)$ satisfying the usual conditions, and 
$\sigma(x)$=$(\sigma^i_j(x))$ is a $C^2$~ $n\times m$ matrix-valued diffusion coefficient function such that the
$C^2$~$n\times n$ matrix-valued function $a(x)=\sigma(x)\sigma^*(x)$ is nondegenerate for all $x$.

Let us assume for simplicity (although one can work without this assumption) that solutions of SDE~\eqref{eq:sde-perturbation} are non-explosive, i.e., for any initial condition $X_\eps(0)=x_0$,
the solution process $X_\eps(t)$ is defined for all times $t\ge0$ with probability~1.
These solutions can be described as 
continuous Markov processes
with generator $A_\eps$ whose action on smooth functions~$f$ with compact support is given by
\begin{align*}
A_\eps f (x)&=\langle b(x),Df(x)\rangle+\frac{\eps^2}{2}\Tr(a(x)D^2f(x))
\\&
=\sum_{i=1}^n b^i(x)\partial_i f(x)+\frac{\eps^2}{2}\sum_{i,k=1}^n a^{ik}(x)\partial_{ik}f(x),\quad x\in\R^n.
\end{align*}
where $Df(x)=(\partial_1f(x),\ldots,\partial_nf(x))$,
$D^2f(x)=(\partial_{ik}f(x))_{i,k=1}^n$, and the angular brackets denote the standard inner product in $\R^n$.
The distribution of process $X_\eps$ started at time $0$ at a point $x$ will be denoted by~$\Pp_{\eps,x}$.

We will observe the process $X_\eps$ only while it evolves within a domain~$\Domain\subset\R^n$. In this section, 
it does not
have to be bounded, and its boundary $\partial \Domain$ is not required to be smooth, except that it is required to contain a $C^2$ 
hypersurface~$M$ such that $\Domain\cup M$ is a path-connected set.

If  $X_\eps(0)=x\in \Domain$, we define $\tau_\eps=\inf\{t\ge 0: X_\eps(t)\in \partial \Domain\}\in(0,\infty]$. 
We are interested in the distribution of the exit point $X_\eps(\tau_\eps)$ conditioned on the
event~$C_\Gamma=\{\tau_\eps<\infty, X_\eps(\tau_\eps)\in \Gamma\}$, where $\Gamma$ is a subset
of $\partial\Domain$ containing~$M$.

The assumptions that we have made guarantee that for any $\eps>0$ and any $x\in\Domain$, under $\Pp_{\eps,x}$
the diffusion $X_\eps$ conditioned on $C_\Gamma$ and stopped at~$\partial\Domain$ is also a 
diffusion process. The generator $A_{\Gamma,\eps}$ of the conditioned diffusion is well-defined on every
$f\in C^2_0(\Domain)$ (which we understand to be defined on the whole $\R^n$) and is given by
\begin{equation}
A_{\Gamma,\eps} f(x)=Lf(x)+\frac{\eps^2}{h^\eps(x)}\sum_{i,j=1}^n a^{ij}(x)\partial_j h^\eps(x)\partial_i f(x),
\quad x\in \Domain,
\label{eq:conditioned-generator-eps}
\end{equation}
where $h^\eps$ defined by
\[
h^\eps(x)=\Pp_{\eps,x}(C_\Gamma),\quad x\in \bar{\Domain},
\]
is $C^2$ and strictly positive in~$\Domain$ for all~$\eps>0$. In other words, the conditioned
diffusion has the diffusion coefficients of the original unconditioned diffusion, but the drift
coefficient $\bar b_\eps$ of the conditioned diffusion is given by 
\begin{equation}
\label{eq:conditioned-drift-eps}
\bar b_\eps(x)=b(x)+\eps^2a(x)\frac{Dh^\eps(x)}{h^\eps(x)},\quad x\in \Domain.
\end{equation}

Formulas~\eqref{eq:conditioned-generator-eps} and~\eqref{eq:conditioned-drift-eps} can be viewed
as specific cases of Doob's $h$-transform. Although they are well-known and 
valid under very mild assumptions, no rigorous
and complete exposition of Doob's $h$-transform for conditioned diffusions is known to us. 
Section~\ref{sec:h-transform} aims to be such an exposition containing some relevant results
that are rigorous and general. Specifically, Lemma~\ref{lem:h-is-C2}
implies that $h^\eps$ is differentiable (so~\eqref{eq:conditioned-generator-eps} and~\eqref{eq:conditioned-drift-eps}
make sense), and Theorem~\ref{thm:h-transform-for-diffusions} implies the rest of the claims made in the last
paragraph.

\medskip

Let us now study the effect of some natural assumptions on the limiting behavior of the drift $\bar b_\eps$ of the conditioned diffusion
introduced in~\eqref{eq:conditioned-drift-eps}.

Let us suppose that there is an open set $G$, a point $x_0\in \Domain\cap G$ and a vector field $ \bar b_0\in C^2(G)$  satisfying the following properties:
\begin{itemize}
\item[(A1)] Let $(S^t)_{t\ge 0}$ denote the flow generated by $\bar b_0$. We assume that there are $T>0$ and
$z\in M\cap G$
 such that
\begin{enumerate}[(i)]
\item $S^tx\in \Domain\cap G$ for all $t\in[0,T)$;
\item $z=S^Tx_0$;
\item $\bar b_0$ is transversal to  $M$
at $z$, i.e., $\bar b_0(z)$ does not belong to the tangent hyperplane $T_zM$.
\end{enumerate}
\item[(A2)] We assume that $\bar b_\eps$ admits a $C^2$ continuation onto $G$ and there are positive constants $C$ and $\alpha$ such that for $\eps\in(0,1)$,
\[
\left|\bar b_\eps(x) - \bar b_0(x)\right|\le C\eps^{1+\alpha},\quad x\in G.
\]
\end{itemize}

The first of these properties is called the Levinson condition. It describes the mutual geometry of
the vector field $\bar b_0$ and the domain~$\Domain$. The second property also involves $\bar b_\eps$ and
states that it converges to $\bar b_0$ uniformly and sufficiently fast.

Assumptions~(A1) and~(A2) are imposed so that we can directly apply
the main result of~\cite{B:MR2739004} that states that if in the Levinson case the deterministic
perturbation of the drift converges to zero faster than the stochastic perturbation
and the perturbation of the initial condition (the latter is identically zero in our situation, although one can easily
consider nonzero perturbations of initial conditions as well), then,  under appropriate rescaling, the exit
location has asymptotically Gaussian distribution. Moreover, an explicit expression for the latter
is available, so let us introduce the relevant notation and state one of our main results. 

Let $I$ be the $n\times n$ identity
matrix. Let $\Phi(t)$ be the linearization
of the flow~$S$ along the orbit of $x_0$:
\[
\dot \Phi(t)=D\bar b_0(S^tx_0)\Phi(t),\quad \Phi(0)=I.
\]
Also, for any vector $v\in\R^n$ we uniquely decompose it into
\[
v=\pi_{b}\cdot \bar b_0(z)+\pi_Mv,
\]
where $\pi_{b} v\in\R$ and $\pi_Mv\in T_zM$.

The following is a direct consequence of Theorem~\ref{thm:h-transform-for-diffusions} on $h$-transform
and the main result of~\cite{B:MR2739004}.

\begin{Theorem}
\label{thm:conditioned-exit-CLT}
Under assumptions~{\rm(A1)} and~{\rm(A2)}, the conditional distribution of
$\eps^{-1}(\tau_\eps-T,X_\eps(\tau_\eps)-z)$ given~$C_\Gamma$ converges weakly (as $\eps\to0$), 
to a nondegenerate $n$-dimensional
Gaussian distribution given by $(-\pi_{b}\phi,\pi_M\phi)$,
where
\[
\phi=\Phi(T)\int_0^T\Phi^{-1}(s)\sigma(S^tx_0)dW(t).
\]
\end{Theorem}

Although this is our central result, it would be useless unless we explain why 
assumptions~(A1) and~(A2) hold true in a large class of cases. 
Of course, if $h^\eps\equiv1$, then the conditioning is trivial
and if the original vector field satisfies the Levinson condition, then Theorem~\ref{thm:conditioned-exit-CLT} applies.
We are interested though in conditioning on unlikely events, and we begin with a simple guiding example.

\section{Example}\label{sec:example}

In this short section we consider a simple guiding example where all calculations can be done explicitly.
In this example, $\Domain=\{x=(x_1,x_2)\in\R^2:\ x_1>0\}$,  $b(x)\equiv(b^1,b^2)$, where $b^1>0$ and
$b^2$ is any constant, and $a(x)$ is the identity matrix for all $x$. In other words, this is Brownian motion in half plane, with drift directed away from the boundary 
of the half-plane and can be described by the following stochastic equation with additive noise:
\begin{align*}
dX^1_\eps(t)&=b^1dt+\eps dW^1(t),\\
dX^2_\eps(t)&=b^2dt+\eps dW^2(t).
\end{align*}
For any initial condition $x\in\Domain$, the process $X_\eps$ reaches $\partial\Domain=\{(x_1,x_2)\in\R^2:\ x_1=0\}$ with positive probability
$h^\eps(x_1,x_2)$. Due to the translational invariance of the system along the $x_2$-axis, this probability depends only on~$x_1$,
$h^\eps(x_1,x_2)=h^\eps(x_1)$, and this function of one variable is $C^2$ and satisfies
\begin{align*}
&b^1 \partial_{1}{h^\eps}(x_1)+\frac{\eps^2}{2}\partial_{11}{h^\eps}(x_1)=0,\quad x_1>0,\\
&h^\eps(0)=1,\\
&h^\eps(x_1)\to 0, \quad x_1\to\infty.
\end{align*}
Solving this linear ODE, we obtain 
\[
h^\eps(x_1,x_2)=e^{-\frac{2b^1 x_1}{\eps^2}},\quad x_1\ge 0.
\]
Applying Doob's $h$-transform (Theorem~\ref{thm:h-transform-for-diffusions}) we obtain that the process $X_\eps$
conditioned to hit $\partial \Domain$ is a diffusion with the same diagonal diffusion matrix and drift given by
\begin{align*}
\bar b_{\eps}^1(x)=&b^1+\eps^2\frac{-\frac{2b_1}{\eps^2}e^{-\frac{2b^1 x_1}{\eps^2}}}{e^{-\frac{2b^1 x_1}{\eps^2}}}=-b^1,\\
\bar b_{\eps}^2(x)=&b^2,\quad x\in\Domain.
\end{align*}
Therefore, $\bar b_\eps\equiv (-b^1,b^2)$ for all $\eps$. The conditioned diffusion is thus
also Brownian motion with drift $(-b^1,b^2)$. This drift is directed towards the boundary and does not depend on $\eps$,
so conditions (A1) and (A2) follow and Theorem~\ref{thm:conditioned-exit-CLT} applies. The Gaussian vector
$\phi$ of Theorem~\ref{thm:conditioned-exit-CLT} can also be easily computed explicitly since $\bar b_0\equiv (-b^1,b^2)$, $D\bar b_0=0$, and 
$\Phi(t)=I$ for all $t$.

In general, when the drift $b$ is not constant and boundaries are curved, we do not expect the conditional drift do be independent
of $\eps$. However, in that case one can use PDE methods to check that conditions (A1) and (A2) still hold.

\section{PDE approach to checking conditions of Theorem~\ref{thm:conditioned-exit-CLT}}
\label{sec:PDE-main-results}
From now on, besides the assumptions on $\Domain$ and $\Gamma$  we have made earlier, we assume that 
the domain $\Domain$ is bounded and $\Gamma$ is open in the relative topology of $\partial\Domain$,
although some of our arguments can be modified to fit other situations as well.

According to Lemma~\ref{lem:Lh}, the function $h^\eps$ satisfies
\begin{equation}
\label{eq:Lh-eps}
\sum_{i=1}^n b^i(x)\partial_i h^\eps(x)+
\frac{\eps^2}{2}\sum_{i,j=1}^n a^{ij}(x)\partial_{ij}h^\eps(x)=0,\quad x\in\Domain.
\end{equation}
We also know that on the boundary
\begin{equation*}
h^\eps(x)=\begin{cases}
           1,& x\in\Gamma,\\
           0,& x\in\partial\Domain\setminus\Gamma,
          \end{cases}
\end{equation*}
but continuity of $h^\eps$ at a boundary point $x$ depends on whether $x$ is {\it regular}, i.e.,
whether for all $\eps,t>0$
\[
 \lim_{\substack{\Domain \ni y\to x}} \Pp_{\eps,x}\{\tau_\eps>t\}=0.
\]
In fact (see, e.g., implication (i)$\Rightarrow$(ii) of Theorem~2.3.3 in~\cite{Pinsky:MR1326606})
\begin{equation}
\label{eq:h-continuous-at-boundary}
 \lim_{\Domain\ni y\to x} h^\eps(y)=\begin{cases}
           1,& \mbox{regular}\ x\in\Gamma,\\
           0,& \mbox{regular}\ x\in\partial\Domain\setminus\bar \Gamma,
          \end{cases}
\end{equation}
so let us make an assumption that all boundary points are regular.
This condition is implied by 
the following {\it external cone condition}: for every boundary point $x$ there are a cone $K$
with base $x$ and a neighborhood $U$ of $x$ such that $K\cap U\cap \Domain=\emptyset$, see, e.g.,~\cite[Section 2.3]{Pinsky:MR1326606}.

Since $0<h^\eps(x)\le 1$ for $x\in \Domain\cup \Gamma$, we can make a Hopf--Cole type logarithmic
change of variables:
\[
v^\eps(x)=-\eps^2\log h^\eps(x),\quad x\in \Domain\cup \Gamma,
\]
so that $0\le v^\eps(x)<\infty$ for $x\in \Domain\cup \Gamma$. 
The reason for using this transformation is that now $\bar b_\eps$ can be represented as
\begin{equation}
\label{eq:conditioned-drift-eps-via-Hopf-Cole}
\bar b_\eps(x)=b(x)-a(x)Dv^\eps(x),\quad x\in \Domain,
\end{equation}
and one may hope to study the behavior of $\bar b_\eps$ via the analysis of $Dv^\eps$.
Plugging $h^\eps=\exp\{-v^{\eps}/\eps^2\}$ into~\eqref{eq:Lh-eps} and \eqref{eq:h-continuous-at-boundary}, we obtain
\begin{equation}
\label{eq:viscous-HJB}
-\frac{\eps^2}{2}\sum_{i,j=1}^n a^{ij}\partial_{ij}v^\eps-\sum_{i=1}^n b^i\partial_i v^\eps+
\frac{1}{2}\sum_{i,j=1}^n a^{ij}\partial_{i}v^\eps\partial_{j}v^\eps=0\quad\mbox{in}\ \Domain,
\end{equation}
and
\begin{equation}
\label{eq:boundary-cond-for-v-eps}
 \lim_{\Domain\ni y\to x} v^\eps(y)= \begin{cases}
0,& x\in \Gamma,
\\ +\infty,& x\in\partial\Domain\setminus\bar \Gamma.
\end{cases}
\end{equation}
Equation~\eqref{eq:viscous-HJB} is a stationary HJB equation with positive viscosity and, since the boundary
condition~\eqref{eq:boundary-cond-for-v-eps} does not depend on $\eps$ one should expect 
that
solutions of such equations converge as $\eps\to 0$ to a solution $v^0$ of the inviscid HJB equation
\begin{equation}
\label{eq:HJB-general}
-\sum_{i=1}^n b^i\partial_i v^0+
\sum_{i,j=1}^n a^{ij}\partial_{i}v^0\partial_{j}v^0=0,
\end{equation}
equipped with boundary conditions
\begin{equation}
\label{eq:boundary-cond-for-HJB}
v^0(x)=\begin{cases}
        0,& x\in\Gamma,\\
        +\infty,& x\in \partial\Domain\setminus\bar\Gamma.
       \end{cases}
\end{equation}

Therefore, a natural candidate for the vector field $\bar b_0$ of Conditions~(A1) and~(A2)
is $b-aDv^0$. 
However, several difficulties arise. First of all, classical smooth solutions of~\eqref{eq:HJB-general}--\eqref{eq:boundary-cond-for-HJB}
often do not exist in the entire domain, and one has to deal with generalized solutions. Viscosity solutions 
(see, e.g.,~\cite{Lions:MR667669,CIL:MR1118699,Fleming-Soner:MR2179357})
form a natural class of solutions, however even then the boundary condition ~\eqref{eq:boundary-cond-for-HJB} cannot be satisfied
on the entire boundary and must be understood in a generalized sense. Still, one can prove convergence $v^\eps\to v^0$, where $v^0$ is the value function for a variational problem associated with the HJB
equation ~\eqref{eq:HJB-general}--\eqref{eq:boundary-cond-for-HJB} and is a viscosity solution of this equation (see \cite{F2:MR512217}, \cite{EI:MR781589}). 
For singular perturbation problems one can
establish estimates on convergence rates (see, e.g.,~\cite[Chapter 6]{Lions:MR667669}) and even establish
uniform asymptotic expansions of the form
\[
 v^\eps=v^0+\eps^2 v_1+\eps^4v_2+\ldots \eps^{2k}v_k+o(\eps^{2k})
\]
in regions of strong regularity of the solution $v^0$, where it is smooth and coincides with the classical solution
obtained by the method of characteristics, see~\cite{F1:MR0304045,FS:MR833404}. However, we need estimates of this  
kind not just for $v^\eps$, but also for~$Dv^\eps$. Such estimates have been obtained in ~\cite{F1:MR0304045} for evolutionary problems with continuous 
boundary values.
Our case is more delicate and no appropriate results seem to exist in the literature. Major difficulties come from the nature of the
boundary condition and the fact that on the level of stochastic control our problem corresponds to an infinite horizon one.

In the remainder of this section we will provide sufficient conditions for estimate
\begin{equation}
\label{eq:expansion-for-Dv-eps}
Dv^\eps=Dv^0+\eps^2 Dv_1+o(\eps^2)  
\end{equation}
to hold uniformly in regions of strong regularity.

From now on we will restrict ourselves to the
case of isotropic additive noise and assume that $a(x)$ is the identity matrix, i.e., $a^{ij}(x)\equiv\delta^{ij}$
(although the choice of $\sigma(x)$ guaranteeing this is not unique, the distribution of the resulting diffusion process does not depend on that choice, so
for definiteness one may choose $\sigma(x)$ to be the $n\times n$ identity matrix).
We will also assume that $b$ is smooth (i.e. $b\in C^\infty(\bar\Domain)$), and that $\Domain$ is a bounded
domain with smooth boundary $\partial\Domain$. These will be our standing assumptions that will not be repeated in the rest of the paper.
The smoothness assumptions may not be necessary. However, we
impose them for simplicity and to be able to cite explicitly results from papers where they were used. We remark
that, requiring less regularity for the region of strong regularity, to carry out the program of our paper it would be enough 
to assume that $b\in C^4(\bar\Domain)$.

Equation~\eqref{eq:viscous-HJB} now becomes
\begin{equation}\label{a2}
-\frac{\eps^2}{2}\Delta v^\eps-\langle b,Dv^\eps\rangle +\frac{1}{2}|Dv^\eps|^2=0,
\end{equation}
where $\Delta$ is the Laplace operator: $\Delta f=\partial_{11}f+\ldots+\partial_{nn}f$, and $|\cdot|$ is the Euclidean norm.
Equation~\eqref{eq:HJB-general} becomes
\begin{equation}\label{a3}
-\langle b,Dv\rangle +\frac{1}{2}|Dv|^2=0.
\end{equation}

It is well-known that natural candidates for viscosity solutions of HJB equations are value functions
of the associated control/variational problems. Let us make it precise for equation~\eqref{a3}. To that end 
we make the following additional assumptions on $b$:
\begin{itemize}
\item[(B1)]
There are a constant $\lambda>0$ and a relatively open set $U\subset \bar\Domain$ such that $\bar\Gamma\subset U$
and
\begin{equation}\label{a4}
\langle b(x),\nu(x)\rangle\leq -\lambda<0,\quad x\in U\cap\partial\Domain,
\end{equation}
 where for a point
$x\in\partial\Domain$, $\nu(x)$ denotes the outward unit normal vector to $\partial\Domain$ at $x$. 
\item[(B2)] 
For any absolutely continuous path  $\gamma(\cdot):[0,+\infty)\to \bar\Domain$,
\[
\int_0^{+\infty}|\dot\gamma(t)-b(\gamma(t))|^2dt=+\infty.
\]
\end{itemize}
The assumption (B2) was introduced in \cite{F2:MR512217} and was also used in \cite{EI:MR781589}.

For a curve $\gamma\in H^1_{\rm loc}([0,+\infty);\R^n)$ we set
\[
\tau_\gamma=\inf\{t\ge 0: \gamma(t)\in \partial \Domain\},
\]
\[
\bar\tau_\gamma=\inf\{t\ge 0: \gamma(t)\in \R^n\setminus\bar\Domain\},
\]
and for $T>0$ and a curve $\gamma\in H^1([0,T];\bar\Domain)$ we define the action functional
\begin{equation}
\label{eq:action}
 A^{0,T}(\gamma)=\frac{1}{2}\int_0^{T}|\dot\gamma(t)-b(\gamma(t))|^2dt,
\end{equation}
and the value function $v^0:\bar\Domain\to\R$ for variational problem associated with the HJB equation
\begin{equation}
\label{eq:functionv0}
v^0(x)=\inf \Bigl\{A^{0,\bar\tau_\gamma}(\gamma):\  \gamma\in H^1_{\rm loc}([0,+\infty);\R^n),\ \gamma(0)=x,\ \gamma(\bar\tau_\gamma)\in\Gamma \Bigr\}.
\end{equation}

\begin{Lemma}\label{lem:existence-of-RSR}
Under assumptions {\rm(B1)} and {\rm(B2)}, there is a relatively open subset  $N$ of
$\Domain\cup\Gamma$  satisfying the following conditions:
\begin{enumerate}[\rm (a)]
\item $\bar N\subset(\Domain\cup\Gamma)\cap U$.
 \item For every $x\in N$, there is a minimizing curve providing minimum in~\eqref{eq:functionv0}. If $\gamma_x^1, \gamma_x^1$ are two minimizing curves then $\tau_{\gamma_x^1}
=\bar\tau_{\gamma_x^1}=\tau_{\gamma_x^2}=\bar\tau_{\gamma_x^2}$ and $\gamma_x^1(t)= \gamma_x^1(t), t\in[0,\tau_{\gamma_x^1}]$. In this sense
the minimizing curve is unique and we will denote it by $\gamma_x$.
\item For every $x\in N$ and all $t\in[0,\tau_{\gamma_x}]$, $\gamma_x(t)\in N$.
\item The function $v^0\in C^\infty(N)$. 
\item
For every $x\in N$,
\[
\dot \gamma_x(t)=-b_0(\gamma_x(t)), \quad t\in [0,\tau_{\gamma_x}], 
\]
where 
\begin{equation}
b_0(y)= -b(y)+Dv^0(y),\quad y\in N,
\label{eq:choice-of-b_0}
\end{equation}
and $b_0(\gamma_x(\tau_{\gamma_x}))$ is transversal to $\Gamma$.
The curve $\varphi_x(t):=\gamma_x(\tau_{\gamma_x}-t), t\in[0,\tau_{\gamma_x}]$,
is {\it the characteristic curve passing through $x$} and it solves the equation $\dot \varphi_x(t)=b_0(\varphi_x(t))$.
\item For every characteristic curve $\varphi(t), 0\leq t\leq s$, in $N$ there exists an open set 
$B\subset\R^{n-1}$ 
and a local smooth parametrization
$\psi:B\to \Gamma$ such that $\psi(y_0)=\varphi(0)$ for some $y_0\in B$ and if $x(t,y)$ denotes the characteristic curve such that $x(0,y)=\psi(y)$
(in particular, $\varphi(t)=x(t,y_0)$),
then $x(t,y)$ is a smooth diffeomorphism of $[0,s]\times B$ onto its image.
\end{enumerate}
\end{Lemma}

Any region $N$ satisfying the conditions (a)--(f) of Lemma~\ref{lem:existence-of-RSR} will be called a {\it region of strong regularity}.
Our definition of a region of strong regularity is designed specifically for our case and just lists all the properties such a region
should possess. In the literature, the term ``region of strong regularity'' may refer to a slightly different and more general object, 
see ~\cite{F1:MR0304045,FS:MR833404}.
We prove Lemma~\ref{lem:existence-of-RSR} stating
the existence of regions of strong regularity in Section~\ref{subsec:assumpprel}. Although the method employed there gives only 
existence in principle, in reality regions of strong regularity may be quite large.

Part (e) of Lemma~\ref{lem:existence-of-RSR} implies that condition~(A1) holds true for any $x_0$ in a region of strong regularity.
Let us now address condition~(A2) and introduce our last assumption.
\begin{itemize}
 \item[(C)] The set $\Gamma$ is flat, i.e., it is contained in a hyperplane, and the entire domain~$\Domain$ lies entirely
 on one side of that hyperplane.
\end{itemize}

\begin{Theorem}\label{thm:asymptotics-for-Dv-eps} Suppose assumptions {\rm(B1)}, {\rm(B2)}, and~{\rm(C)}
hold. Let $N$ be a region of strong regularity. 
Then there is a function $v_1\in C^\infty(N)$ such that~\eqref{eq:expansion-for-Dv-eps}
holds uniformly on compact subsets of $N$.
\end{Theorem}

Section~\ref{sec:PDE} is devoted to the proof of Theorem~\ref{thm:asymptotics-for-Dv-eps}. This is the most technical
part of this paper. Our multi-stage hard-analysis proof is based on
the approach developed in~\cite{F1:MR0304045} and~\cite{FS:MR833404}, and involves PDE and stochastic methods.

The following is a simple corollary of Lemma~\ref{lem:existence-of-RSR} and Theorem~\ref{thm:asymptotics-for-Dv-eps}.
\begin{Theorem} Suppose assumptions {\rm(B1)}, {\rm(B2)}, and~{\rm(C)} hold.
We recall that $\bar b_\eps=b-Dv^\eps$, see~\eqref{eq:conditioned-drift-eps-via-Hopf-Cole}. Also, let
$\bar b_0=-b_0$, where $b_0$ is defined in~\eqref{eq:choice-of-b_0}.
Then conditions~{\rm(A1)} and~{\rm(A2)} hold for any point $x_0$
belonging to a region of strong regularity and vector fields $\bar b_\eps$ and $\bar b_0$.
 In particular, the CLT for conditioned exit distributions 
stated in Theorem~\ref{thm:conditioned-exit-CLT} applies.
\end{Theorem}

While assumptions (B1) and (B2) define the setting,  the technical assumption (C) does not seem to be
necessary, and we believe that one can obtain similar results without it. However, our current approach depends on the additive
noise character of equation ~\eqref{a2} and does not work for curved boundaries since flattening transformations do not preserve its additive noise structure.

\section{Asymptotics for the elliptic HJB equation}\label{sec:PDE}

The goal of this section is to prove Lemma~\ref{lem:existence-of-RSR} and Theorem~\ref{thm:asymptotics-for-Dv-eps}.
\subsection{Regions of strong regularity. Proof of Lemma~\ref{lem:existence-of-RSR}} 
\label{subsec:assumpprel}
We recall the following results from \cite{EI:MR781589,F2:MR512217}.
\begin{Theorem}\label{thm:1}
Let 
{\rm (B2)} be satisfied. 
Then:
\begin{itemize}
\item[(i)] For every region $\Domain'$ such that $\bar{\Domain'}\subset\Domain\cup\Gamma$ there exists
a constant $C(\Domain')$ independent of $\eps$ such that
\begin{equation}
\label{a5}
\sup_{x\in\Domain'}(|v^\eps(x)|+|Dv^\eps(x)|)\leq C(\Domain').
\end{equation}
\item[(ii)] The function $v^0$ defined by ~\eqref{eq:functionv0} is Lipschitz on $\bar\Domain$ and
\[
 v^0(x)=\inf \Bigl\{A^{0,\tau_\gamma}(\gamma):\  \gamma\in H^1_{\rm loc}([0,+\infty);\R^n),\ \gamma(0)=x,\ \gamma(\tau_\gamma)\in\Gamma \Bigr\},
\]
on $\Domain\cup\Gamma$. It is
 a viscosity solution of
\begin{equation}\label{a3a}
\begin{cases}
-\langle b,Dv^0\rangle +\frac{1}{2}|Dv^0|^2=0\quad\mbox{in}\,\,\Domain,
\\
v^0=0\quad\mbox{on}\,\,\Gamma.
\end{cases}
\end{equation}
\item[(iii)]Uniformly on compact subsets of $\Domain\cup\Gamma$, 
\begin{equation*}
\lim_{\eps\to 0}v^\eps=v^0.
\end{equation*}
\end{itemize}
\end{Theorem}
Part (i) is proved in Lemma~2.2 of~\cite{EI:MR781589} and part~(ii) in Lemmas~2.3 and 2.4 of~\cite{EI:MR781589}. 
Part~(iii) was first proved in~\cite{F2:MR512217} and later a more PDE proof was given in Theorem~2.1 of ~\cite{EI:MR781589}. Condition~(B1) 
is not needed in Theorem~\ref{thm:1}. 

We denote the Euclidean distance in $\R^n$ and the induced Hausdorff distance on subsets of $\R^n$  by $\dist(\cdot,\cdot)$.
We denote
\[
\Domain_\delta=\{x\in\Domain:\dist(x,\partial\Domain)<\delta\},\quad d(x)=\dist(x,\partial\Domain).
\]
We recall that $d\in C^1(\Domain_{\delta})$ for some $\delta>0$, and $Dd(x)=-\nu(x)$ on $\partial\Domain$. Thus
there exist $c,\alpha,\delta>0$ such that the function $\tilde d(x)=cd(x)$ satisfies
\begin{equation}
\label{a6}
-\langle b(x),D\tilde d(x)\rangle +\frac{1}{2}|D\tilde d(x)|^2\leq-\alpha<0,\quad
x\in U\cap \Domain_{\delta}.
\end{equation}

The proof of Lemma~\ref{lem:existence-of-RSR} is based on the method of characteristics. We refer the reader to
\cite[Chapter~3]{Evans:MR2597943},
\cite{Caratheodory-I:MR0192372},
\cite{Garabedian:MR1657375},
\cite{Maurin:MR0507436}, and Appendices in \cite{FS:MR833404} and \cite{F1:MR0304045} for the overview of the method of
characteristics. We sketch the proof below.

Let $\psi:B\to\Gamma$ be a local smooth parametrization of $\Gamma$, where $B\subset\R^{n-1}$ is an open set, and let $x_0=\psi(y_0)\in \psi(B)$.
It is easy to see that any solution of
~\eqref{a3} must satisfy $Du(\psi(y)))=
2\langle b(\psi(y)),\nu(\psi(y))\rangle\nu(\psi(y))\not= 0$. The characteristic system for ~\eqref{a3} 
is the following (see \cite[Chapter~3]{Evans:MR2597943}). Functions
\[
 x(t):=x(t,y),\quad  z(t):=z(t,y), \quad  p(t):=p(t,y)
\]
satisfy
\begin{equation}\label{charactsystem}
\begin{cases}
\dot p(t)=\langle (Db(x(t))^*,p(t)\rangle,
\\
\dot z(t)=\langle -b(x(t))+p(t),p(t)\rangle,
\\
\dot x(t)= -b(x(t))+p(t),
\end{cases}
\end{equation}
with initial conditions
\[
 x(0)=\psi(y),\quad z(0)=0,\quad p(0)=2\langle b(\psi(y)),\nu(\psi(y))\rangle\nu(\psi(y)).
\]
The curve $x(t)$ is called the characteristic starting at $\psi(y)$. Since the right hand side of
~\eqref{charactsystem} is a smooth function of $x,z,p$, the system has a unique local solution for any initial conditions
which is smooth with respect to $t$ and the initial conditions. Since $b$ and $\psi$ are smooth, we thus get that
$x(t,y),z(t,y),p(t,y)$ are smooth functions of $t,y$ in some neighborhood of
$(0,y_0)$. Since
\[
 \langle \dot x(0,y),\nu(\psi(y))\rangle=\langle b(\psi(y)),\nu(\psi(y))\rangle\not = 0,
\]
$\dot x(0,y)$ is transversal to $\Gamma$ at every point $\psi(y)$ (the so called non-characteristic condition), and thus ${\rm det}(Dx(0,y_0))\not = 0$,
which implies that $x(t,y)$ is a local smooth diffeomorphism, and hence the inverse functions $t(x),y(x)$ 
exist and are smooth. Therefore $u(x):=z(t(x),y(x))$ is a smooth function in some relatively open
set containing $x_0$.  One then proves (see \cite[Chapter~3]{Evans:MR2597943})
that $Du(x)=p(t(x),y(x))$ and $u$ satisfies ~\eqref{a3} in this set. Thus for every $x_0\in \Gamma$, $u$ is $C^\infty$ in some relatively open region
$C_{\tilde B,\psi,\tau}=\{x(t,y):y\in\tilde B,0\leq t<\tau\}$ for some $\tau>0$ and an open set $\tilde B\subset\R^{n-1}$. The union
of such sets for $x_0\in \Gamma$ will give us a relatively open region $\tilde N$ of $\bar\Domain$ such that $\tilde N\cap \partial\Domain=\Gamma$, and
$u\in C^\infty(\tilde N)$ and satisfies ~\eqref{a3}.
We now claim that
\begin{equation*}
\label{a7}
u(x)=\min_{\mathcal{U}_{x,\tilde N}}A^{0,\tau_\gamma}(\gamma),
\end{equation*}
where
\begin{multline*}
\mathcal{U}_{x,\tilde N}=\{\gamma(\cdot)\in H^1_{\rm loc}([0,+\infty);\mathbb{R}^n):\\ \gamma(0)=x,\ \gamma(\tau_\gamma)\in\Gamma,\ \tau_\gamma<+\infty,\ 
\gamma(s)\in\tilde N,\  0\leq s\leq \tau_\gamma\},
\end{multline*}
and the minimizing curve is unique in $\mathcal{U}_{x,\tilde N}$. We will use that for $p\in\R^n$,
\[
 \frac{1}{2}|p|^2=\sup_{y\in\R^n}\{\langle y,p\rangle-\frac{1}{2}|y|^2\}
\]
with equality if and only if $y=p$.
Let $\gamma(\cdot)\in \mathcal{U}_{x,\tilde N}$. Then
\[
 \begin{split}
 &
u(x)=u(x)-u(\gamma(\tau_\gamma))=-\int_0^{\tau_\gamma}\langle Du(\gamma(t)),\dot\gamma(t)\rangle dt
\\
&
=\int_0^{\tau_\gamma}\left(-\langle b(\gamma(t)),Du(\gamma(t))\rangle +\langle Du(\gamma(t)), b(\gamma(t))-\dot\gamma(t)\rangle 
-\frac{1}{2}| b(\gamma(t))-\dot\gamma(t)|^2\right)dt
\\
&
+\int_0^{\tau_\gamma} \frac{1}{2}| b(\gamma(t))-\dot\gamma(t)|^2dt\leq \int_0^{\tau_\gamma} \frac{1}{2}| b(\gamma(t))-\dot\gamma(t)|^2dt
 \end{split}
\]
with equality if and only if $b(\gamma(t))-\dot\gamma(t)=Du(\gamma(t))$. Thus the minimizing curve exists and is unique
and thus conditions (c)--(f) are satisfied if $v^0$ there is replaced by $u$ and $N=\tilde N$. 
Moreover, we can assume that 
$\tilde N\subset U\cap \Domain_{\delta}$. 

We now define
\begin{equation*}
\mathcal{U}_{x}=\{\gamma(\cdot)\in H^1_{\rm loc}([0,+\infty);\mathbb{R}^n):\gamma(0)=x,\gamma(\tau_\gamma)\in\Gamma,\tau_x<+\infty\}.
\end{equation*}
We will show that $u=v^0$ in a subset of $\tilde N$.

\begin{Lemma}
\label{lem:1}For every $\sigma>0$
 there is $\delta(\sigma)>0$ such that
$u=v^0$ in $N_\sigma=\{x\in \tilde N\cap\Domain_{\delta(\sigma)}:{\rm dist}(x,\partial\tilde N\setminus\partial\Domain)
>\sigma\}$.
\end{Lemma}
\bpf
Let $x\in N_\sigma$. It is enough to show that if $\gamma(\cdot)\in \mathcal{U}_{x}\setminus\mathcal{U}_{x,\tilde N}$, then
\begin{equation}
\label{a10}
\frac{1}{2}\int_0^{\tau_\gamma}|\dot\gamma(s)-b(\gamma(s))|^2ds>\tilde u(x).
\end{equation}
Let $t\in(0,\tau_\gamma)$ be the smallest $t$ such that $\gamma(t)\in \Domain\setminus\tilde N$. Then, using (\ref{a6}),
\[
\begin{split}
&
\tilde d(\gamma(t))-\tilde d(x)=\int_0^t\langle D\tilde d(\gamma(s)),\dot\gamma(s)\rangle ds
\\
&
=\int_0^t[\langle D\tilde d(\gamma(s)),\dot\gamma(s)-b(\gamma(s))\rangle
+\frac{1}{2}|\dot\gamma(s)-b(\gamma(s))|^2
\\
&
\quad\quad\quad\quad
+\langle D\tilde d(\gamma(s)),b(\gamma(s))\rangle-\frac{1}{2}|\dot\gamma(s)-b(\gamma(s))|^2]ds
\\
&
\geq
\int_0^t[\langle D\tilde d(\gamma(s)),b(\gamma(s))\rangle-\frac{1}{2}|D\tilde d(\gamma(s))|^2]ds-\frac{1}{2}\int_0^t|\dot\gamma(s)-b(\gamma(s))|^2ds
\\
&
\quad\quad\quad\quad
\geq\alpha t-\frac{1}{2}\int_0^t|\dot\gamma(s)-b(\gamma(s))|^2ds.
\end{split}
\]
Since $\tilde d(x)\geq 0$, we obtain
\begin{equation}
\label{a11}
\frac{1}{2}\int_0^{t}|\dot \gamma(s)-b(\gamma(s))|^2ds\geq \alpha t-\tilde d(\gamma(t)).
\end{equation}
On the other hand,
\[
\begin{split}
&
\sigma\leq |\gamma(t)-x|=\left|\int_0^t\dot \gamma(s)ds\right|\leq\int_0^t|\dot\gamma(s)-b(\gamma(s))|ds+\int_0^t|b(\gamma(s))|ds
\\
&
\leq\sqrt t\left(\int_0^t|\dot\gamma(s)-b(\gamma(s))|^2ds\right)^{\frac{1}{2}}+Ct,
\end{split}
\]
where $C>0$ is a constant.
Thus
\[
\frac{1}{2}\int_0^t|\dot\gamma(s)-b(\gamma(s))|^2ds\geq \frac{{(\sigma-Ct)_+}^2}{2t},
\]
where for $a\in\R$, $a_+=\max(a,0)$.
Since there is a constant $C_1>0$ such that 
\begin{equation}
\label{eq:u(x)_le_c_delta}
\tilde u(x)\leq C_1\delta,\quad x\in \tilde N\cap \Domain_\delta,
\end{equation}
we conclude that there are $t_0(\sigma)>0$ and $\delta_0>0$ such that $\delta<\delta_0$ and $t<t_0(\sigma)$ imply
(independently of $x$) inequality~\eqref{a10}.
Hence from now on we can assume that $\delta<\delta_0$ and 
$t_0(\sigma)\le t <\tau_\gamma$.

Let us denote $\eta=\tilde d(\gamma(t))$. Let now $t_1\in[t,\tau_\gamma)$ be such that $\gamma(t_1)\in\partial(\Domain_\eta\cap U)\setminus\partial\Domain$, and
$\gamma(s)\in \Domain_\eta\cap U$ for $t_1<s\leq\tau_\gamma$. Arguing as in the derivation of\eqref{a11}, we obtain
\begin{equation}
\label{a12}
\frac{1}{2}\int_{t_1}^{\tau_\gamma}|\dot\gamma(s)-b(\gamma(s))|^2ds\geq 
\tilde d(\gamma(t_1))+\alpha(\tau_\gamma-t_1).
\end{equation}
Denote $r=\dist(\partial U,\Gamma)$. If $\tilde d(\gamma(t_1))=\eta=\tilde d(\gamma(t))$, then, 
combining~(\ref{a11}) and~(\ref{a12})
we obtain 
\[
\frac{1}{2}\int_{0}^{\tau_\gamma}|\dot\gamma(s)-b(\gamma(s))|^2ds\geq \alpha(t_0(\sigma)+\tau_\gamma-t_1),
\]
and, due to~\eqref{eq:u(x)_le_c_delta}, there is $\delta_1(\sigma)>0$ such that (\ref{a10}) holds true for $\delta<\delta_1(\sigma)$.
If $\tilde d(\gamma(t_1))\ne \eta$, we must have $\gamma(t_1)\in \partial U$, and then
\[
\frac{1}{2}\int_{t_1}^{\tau_\gamma}|\dot\gamma(s)-b(\gamma(s))|^2ds\geq \frac{{(r-C(\tau_\gamma-t_1))_+}^2}{2(\tau_\gamma-t_1)},
\]
which together with~\eqref{eq:u(x)_le_c_delta} implies that there is $\delta_2>0$ such that if $\delta<\delta_2$, 
then~(\ref{a10}) holds unless $(\tau_\gamma-t_1)\geq \tilde t_0>0$ for some $\tilde t_0=\tilde t_0(r)$. In this case,
Combining~(\ref{a11}),(\ref{a12}), and the fact that $\eta\le \delta$,
we obtain that there is $\delta_3(r)>0$ such that for $\delta<\delta_3(r)$,
\[
\frac{1}{2}\int_{0}^{\tau_\gamma}|\dot\gamma(s)-b(\gamma(s))|^2ds\geq \alpha(t_0(\sigma)+\tilde t_0(r))-\eta>\alpha t_0(\sigma).
\]
Using~\eqref{eq:u(x)_le_c_delta} once again
we conclude that there is $\delta_4(\sigma,r)>0$ such that in this case~(\ref{a10}) holds for $\delta<\delta_4(\sigma,r)$. 
The lemma now
follows with $\delta(\sigma)=\min\{\delta_i,\ i=0,\ldots, 4\}$.
\epf

The region $N_\sigma$ itself may not be a region of strong regularity since it may not satisfy condition (f). However we can now
take $N$ to be the union over all $x_0\in \Gamma$ of sets of the form $C_{\tilde B,\psi,\tau}\subset N_\sigma$ described before.

\subsection{Convergence of derivatives}

The main goal of this section is to prove Lemma~\ref{lem:convofderivatives} stating uniform convergence of $Dv^\eps$ to
$Dv^0$ as $\eps\to 0$. Since $b$ and $\partial\Domain$ are smooth it is well known that in this case the function $h^\eps\in C^2(\Omega)$, 
it is continuous at every point of $\Gamma$, $0<h^\eps<1$ in $\Domain$ and, moreover 
(see, e.g.,~\cite[Section 6.5]{Gilbarg-Trudinger:MR1814364},
 $h^\eps$ is smooth in every 
relatively open region
$\Domain'$
such that $\bar{\Domain'}\subset\Domain\cup\Gamma$. Therefore $v^\eps$ is also smooth in every 
relatively open region
$\Domain'$
such that $\bar{\Domain'}\subset\Domain\cup\Gamma$.

Without further loss of generality we suppose that
$\Gamma\subset \{x_n=0\}$ and $N\subset\{x_n\geq 0\}$. Under this assumption, (\ref{a4}) implies that
$-\langle b,e_n\rangle\leq -\lambda<0$  in a relative neighborhood of $\Gamma$ in $\partial \Domain$, where $e_n=(0,...,0,1)$ is the standard basis vector.

For $r>0$, we denote 
\[
N_r=\left\{x=(x',x_n)\in N:{\rm dist}((x',0),\partial N\setminus\Gamma)>2r, x_n<r\right\}.
\] 

\begin{Lemma}
\label{lem:2}
For every $r>0$ there are positive numbers $L,C,\eps_0$ such that for $\eps<\eps_0$,
\begin{equation}
\label{c1}
v^\eps(x)\leq (1+L\eps^2)v^0(x)+Cx_n^2,\quad x\in N_r.
\end{equation}
\end{Lemma}
\bpf
For $y\in N_r$, we denote 
$W_r(y)=\{x=(x',x_n): |x'-y'|<r, x_n<r\}$. 
Let us define
\[
\psi^\eps_y(x)=(1+L\eps^2)v^0(x)+Cx_n^2+C_1x_n|x'-y'|^2,\quad x,y\in N_r.
\]
and show that there is a choice of constants $L,C,C_1$ such that
\[
v^\eps(x)\leq\psi_y^\eps(x),\quad y\in N_r, x\in W_r(y), 
\]
This will obviously imply (\ref{c1}) since for any $x\in N_r$
one can choose $y=x$.

To simplify notation, without loss of generality we will assume that $y'=0$
and use $\psi^\eps$ for $\psi^\eps_0$.

Directly from~\eqref{a3} we know that
\begin{equation}
\label{eq:Dv0-on-Gamma}
Dv^0(x)=2\langle b,e_n\rangle e_n,\quad x\in\Gamma.
\end{equation}
Therefore, we can choose $r$ such that for some $c_0>0$,
\begin{equation}
\label{c2}
|Dv^0(x)|\geq c_0,\quad
\langle -b(x)+Dv^0(x),e_n\rangle\geq c_0, \quad x\in W_r(y),
\end{equation}
and from this point on we write $W$ for $W_r(0)$.
Due to (\ref{a5}), inequality
\begin{equation*}
v^\eps(x)\leq C_1r^2x_n\quad\mbox{on}\,\,\{x\in\partial W:|x'|=r\},
\end{equation*}
 holds with $C_1=C(N)/r^2$.
Since $v^\eps$ converges uniformly to $v^0$ in $N$, the estimate
\begin{equation*}
v^\eps(x)\leq v^0(x)+Cr^2\le \psi^\eps(x)\quad\mbox{on}\,\,\{x\in\partial W:x_n=r\}
\end{equation*}
holds true if we choose, say, $C\geq 1$ in the definition of $\psi^{\eps}$ and let
$\eps$ be small enough. We have
\[
D\psi^\eps(x)=(1+L\eps^2)Dv^0(x)+(2Cx_n+C_1|x'|^2)e_n+2C_1x_n(x',0).
\]
If $L\eps^2\leq 1$ then, due to the smoothness of $v^0$, there is a number $M(C,C_1)>0$ such that
for all $\eps$,
\begin{equation}
\label{c5}
|D\psi^\eps(x)|,\|D^2\psi^\eps(x)\|\leq M(C,C_1),\quad x\in W.
\end{equation}
Here, for an $n\times n$ matrix $A$,
$\|A\|=\max_{i,j=1,\ldots,n}|A_{ij}|$.
 
Using (\ref{c2}) and (\ref{c5}), we obtain
\[
\begin{split}
&-\frac{\eps^2}{2}\Delta \psi^\eps(x)-\langle b(x),D\psi^\eps(x)\rangle +\frac{1}{2}|D\psi^\eps(x)|^2
\geq -\frac{\eps^2nM}{2}
\\
&
+(1+L\eps^2)\left(-\langle b(x),Dv^0(x)\rangle +\frac{1}{2}|Dv^0(x)|^2\right)
+\frac{L\eps^2(1+L\eps^2)}{2}|Dv^0(x)|^2
\\
&
+\langle -b(x)+Dv^0(x),e_n\rangle(2Cx_n+C_1|x'|^2)
+L\eps^2\langle Dv^0(x),e_n\rangle(2Cx_n+C_1|x'|^2)
\\
&
-2C_1x_n\langle b(x),(x',0)\rangle
+2C_1(1+L\eps^2)x_n\langle Dv^0(x),(x',0)\rangle
\\
&
+
\frac{1}{2}\left|(2Cx_n+C_1|x'|^2)e_n+2C_1x_n(x',0)\right|^2
\\
&
\geq
-\frac{\eps^2nM}{2}+c_0(2Cx_n+C_1|x'|^2)+\frac{L\eps^2c_0^2}{2}-C_2x_n|x'|,
\end{split}
\]
where $C_2$ depends only on $C_1$ and the bounds on $|b|,|Dv^0|$. It is now clear that the last line of the above inequality
can be made nonnegative on $W$ by taking $L$ and $C$ large enough. Therefore $\psi^\eps$ is a supersolution
of (\ref{a2}) in $W$ and $v^\eps\leq \psi^\eps$ on $\partial W$ by (\ref{a2}) and (\ref{a3}). Therefore, by
comparison principle (see for instance \cite[Theorem 17.1]{Gilbarg-Trudinger:MR1814364}) 
we obtain $v^\eps\leq\psi^\eps$ in $W$.
\epf

\begin{Corollary}
\label{cor:1}
For every $r>0$ there exist $L_1,\eps_0>0$ such that for $\eps<\eps_0$,
\begin{equation}
\label{c6}
\partial_n v^\eps(x)\leq (1+L\eps^2)\partial_n v^0(x)
\leq  \partial_n v^0(x)+L_1\eps^2,\quad x\in N(r)\cap \Gamma.
\end{equation}

\end{Corollary}
\begin{Corollary}
\label{cor:2}
Let $\eps_0$ be as in Corollary~\ref{cor:1}. Then there is a constant $L_2>0$ such that for $\eps<\eps_0$,
\begin{equation}
\label{c7}
\partial_{nn} v^\eps(x)\leq L_2,\quad x\in N(r)\cap \Gamma.
\end{equation}
\end{Corollary}
\bpf
Let $x\in N(r)\cap \Gamma$. Then
\[
\partial_{nn} v^\eps(x)=\Delta v^\eps(x)
=\frac{2}{\eps^2}\left(-\langle b(x),Dv^\eps(x)\rangle +\frac{1}{2}|Dv^\eps(x)|^2\right).
\]
Since $v_\eps\equiv 0$ on $\Gamma$, we have $Dv^\eps(x)=t(x)e_n$ for some $t(x)\geq 0$. Combining this
with Corollary~\ref{cor:2} and~\eqref{eq:Dv0-on-Gamma}, we obtain 
$t(x)\leq 2\langle b(x),e_n\rangle+L_1\eps^2$.
Hence
\[
\partial_{nn} v^\eps(x)=\frac{2}{\eps^2}\left(-\langle b(x),e_n\rangle t
+\frac{1}{2}t^2\right)\leq \frac{2}{\eps^2}\left(\langle b(x),e_n\rangle L_1\eps^2+\frac{2}{2}L^2\eps^4\right),
\]
which implies our claim.
\epf
\begin{Remark}\rm
\label{rem:1}
A version of Lemma \ref{lem:2} could be obtained without the assumption that $\Gamma$ is flat. Corollaries
\ref{cor:1} and \ref{cor:2} are also true in this case if  $\partial_n v^0$, $\partial_n v^\eps$, and
$\partial_{nn} v^\eps$  are replaced by the respective derivatives along the inward normal vector.
Lemma~\ref{lem:semiconcavity} below is the only place where we have to assume that $\Gamma$ is flat.
\end{Remark}

The following lemma and its proof is a standard semiconcavity result from~\cite{Lions:MR667669} which was also used in \cite{FS:MR833404}.
The difference between it and Lemma~3.1 of~\cite{FS:MR833404} is that we prove the semiconcavity up to 
the boundary~$\Gamma$, although only along coordinate
directions. This is made possible by Corollary~\ref{cor:2}. On compact subsets of $\Domain$ a semiconcavity estimate (\ref{c8}) holds along every direction but we do not need
it here. We include the proof of the lemma for completeness.
\begin{Lemma}
\label{lem:semiconcavity}
If $N'$ is a relatively compact subset of $\Domain\cup\Gamma$, then there is a constant $C(N')>0$ such that
\begin{equation}
 \label{c8}
\partial_{ii}v^\eps(x)\leq C(N'),\quad x\in N',i=1,...,n.
\end{equation}
\end{Lemma}
\bpf
 Since ${\rm dist}(N',\partial\Domain\setminus \Gamma)>0$, there exists a function
$\xi\in C^2_0(\Domain\cup \Gamma)$ 
such that $0\le \xi\le 1$ for all $x$,  $\xi=1$ on $N'$,  and
\begin{equation*}
\frac{|D\xi|^2}{\xi}\leq C_1\quad\mbox{on the support of}\,\,\xi
\end{equation*}
for some $C_1$. To simplify notation, for a function $f$ with real or vector values, we will write $f_{i}$ and $f_{ii}$ for $\partial_i f$ and $\partial_{ii} f$.
Differentiating (\ref{a2}) with respect to $x_i$ we have
\[
-\frac{\eps^2}{2}\Delta v^\eps_i-\langle b_i,Dv^\eps\rangle 
-\langle b,Dv^\eps_i\rangle
+\langle Dv^\eps,Dv^\eps_i\rangle=0,  
\]
\[
-\frac{\eps^2}{2}\Delta v^\eps_{ii}-\langle b_{ii},Dv^\eps\rangle 
-2\langle b_i,Dv^\eps_i\rangle-\langle b,Dv^\eps_{ii}\rangle
+|Dv^\eps_i|^2+\langle Dv^\eps,Dv^\eps_{ii}\rangle=0. 
\]
Therefore, setting $w=\xi v^\eps_{ii}$, we have
\begin{equation}
\label{c10}
 \begin{split}
 &
-\frac{\eps^2}{2}\Delta w+{\eps^2}{\rm tr}\left(\frac{D\xi}{\xi}\otimes 
Dw\right)
=-\frac{\eps^2}{2}\xi\Delta v^\eps_{ii}
-\frac{\eps^2}{2}v^\eps_{ii}\Delta \xi
+{\eps^2}\frac{|D\xi|^2}{\xi}v^\eps_{ii}
\\
&
=\xi\langle b_{ii},Dv^\eps\rangle 
+2\xi\langle b_i,Dv^\eps_i\rangle
+\xi\langle b,Dv^\eps_{ii}\rangle
-\xi|Dv^\eps_i|^2-\xi\langle Dv^\eps,Dv^\eps_{ii}\rangle
\\
&
-\frac{\eps^2}{2}\left(\Delta \xi
-2\frac{|D\xi|^2}{\xi}\right)v^\eps_{ii}
=\xi\langle b_{ii},Dv^\eps\rangle
+2\xi\langle b_i,Dv^\eps_i\rangle
-\xi|Dv^\eps_i|^2
\\
&
+\langle b- Dv^\eps,Dw\rangle
-\langle b- Dv^\eps,D\xi\rangle v^\eps_{ii}
-\frac{\eps^2}{2}\left(\Delta \xi
-2\frac{|D\xi|^2}{\xi}\right)v^\eps_{ii}.
 \end{split}
\end{equation}
The value of $w$ on $\partial\Domain$ is equal to $0$ if $i=1,...,n-1$, and 
it is bounded by Corollary \ref{cor:2} if $i=n$. 
If $w$ has a positive
maximum at a point $x\in\Domain$, then $Dw(x)=0$ and $\Delta w(x)\leq 0$. Since
$|Dv^\eps|$ are bounded on the support of $\xi$, we thus obtain from
(\ref{c10}) that for some numbers $C_2,C_3,C_4$ depending on $\xi$ but not on $\eps$,
\[
(\xi(x))^2(|Dv^\eps_i(x)|^2-C_2|Dv^\eps_i(x)|)\leq C_3+C_4w(x).
\]
Applying the standard inequality $ab<(a^2/c+cb^2)/2$ to $C_2|Dv^\eps_i(x)|$, we obtain
that there are numbers $C_5,C_6>0$ independent of $\eps$ such that
\[
|Dv^\eps_i(x)|^2-C_2|Dv^\eps_i(x)|\ge C_5|Dv^\eps_i(x)|-C_6.
\]
Since $\xi^2$ is bounded by $1$, we obtain
\[
 (w(x))^2\leq \xi(x)^2|Dv^\eps_i(x)|^2\leq (C_6+C_3+C_4w(x))/C_5,
\]
and the resulting uniform upper bound on $w$  implies~\eqref{c8}.
\epf
 
The following lemma is now a consequence of Corollary~\ref{cor:1} and Lemma~\ref{lem:semiconcavity}
(see the proof of Theorem~3.2 of~\cite{FS:MR833404} for a similar argument).
\begin{Lemma}
\label{lem:convofderivatives}
$Dv^\eps$ converges to $Dv^0$ uniformly on every compact subset of~$N$.
\end{Lemma}

\bpf
Let $N'$ be a compact subset of $N$. Since $Dv^0$ is continuous in $N$, it is enough
to show for all $i=1,...,n$  that if for all $m\in\N$, $x^{(m)}\in N'$, and $x^{(m)}\to x\in N'$, $\eps_m\to 0$ as $m\to\infty$, then 
it is possible to extract a subsequence of $(\partial_i v^{\eps_m}(x^{(m)}))$
convergent to $\partial_i v^0(x)$.

 For $h\in \mathbb{R}$
such that $x^{(m)}+he_i\in N'$, (\ref{c7}) implies
\begin{equation}\label{convder1}
v^{\eps_m}(x^{(m)}+he_i)\leq v^{\eps_m}(x^{(m)})+\partial_i v^{\eps_m}(x^{(m)})h+\frac{L_2}{2}h^2.
\end{equation}
Since $\partial_i v^{\eps_m}(x^{(m)})$ is bounded, taking a subsequence if necessary, we have $\partial_i v^{\eps_m}(x^{(m)})\to p$
for some $p$. Since
$v^{\eps_m}\to v^0$ uniformly on $N$,  (\ref{convder1}) implies
\begin{equation}\label{convder2}
v^0(x+he_i)\leq v^0(x)+ph+\frac{L_2}{2}h^2
\end{equation}
for sufficiently small $|h|$ if $x\notin\Gamma$ and for sufficiently small 
$h>0$ if $x\in \Gamma$.
If $i=1,...,n-1$, (\ref{convder2}) clearly implies $p=\partial_i v^0(x)$. If $i=n$ and 
$x\not\in \Gamma$,
we also have $p=\partial_n v^0(x)$. If $x=(x',0)\in\Gamma$, let $x^{(m)}=({x^{(m)}}',a_m)$,
i.e., $a_m=x^{(m)}_n$.
Then $a_m\geq 0$,
$a_m\to 0$ as $m\to\infty$. Using (\ref{c7}), we get
\[
\begin{split}
 &
\partial_n v^{\eps_m}({x^{(m)}}',a_m)-\partial_n v^{\eps_m}({x^{(m)}}',0)
=\int_0^1\frac{d}{dt}\left[\partial_n v^{\eps_m}({x^{(m)}}',ta_m)\right]dt
\\
&
=a_m\int_0^1 \partial_{nn} v^{\eps_m}({x^{(m)}}',ta_m)dt\leq a_mL_2.
\end{split}
\]
Therefore, by (\ref{c6}),
\[
\partial_n v^{\eps_m}(x^{(m)})\leq \partial_n v^{\eps_m}({x^{(m)}}',0)+a_m L_2
\leq 
\partial_n v^0({x^{(m)}}',0)+\eps_m^2L_1+a_m L_2.
\]
This implies $p\leq \partial_n v^0(x)$. Combining this with (\ref{convder2}), we obtain
$p=\partial_n v^0(x)$.
 \epf

\subsection{Asymptotics of $v^\eps$}

We 
follow the method of \cite{FS:MR833404}. We recall that any characteristic $\gamma$
of~\eqref{a3} satisfies $\dot\gamma=b_0(\gamma)$, where
$b_0(x)=-b(x)+Dv^0(x)$.

By the construction in \cite{FS:MR833404}, page 439, for every
$\bar x\in N$ there exists a relatively open in $N$ subregion of regularity
$N_\gamma$ which is a neighborhood of the characteristic curve $\gamma$ connecting
$\bar x$ with $\partial\Domain$ (and consisting of characteristic curves), and
smooth functions $F,G$ on $N_\gamma$ such that
\begin{equation*}
F>0\,\,\mbox{in}\,\,N_\gamma,\quad F=0\,\,\mbox{on}\,\,\partial N_\gamma
\setminus\partial\Domain,
\end{equation*}
\begin{equation*}
G>0\,\,\mbox{in}\,\,N_\gamma\setminus\partial\Domain,\quad G=0\,\,\mbox{on}\,\,
N_\gamma\cap\partial\Domain,
\end{equation*}
\begin{equation*}
-\langle b_0,DF\rangle =1\,\,\mbox{in}\,\,N_\gamma,
\end{equation*}
\begin{equation*}
\langle b_0,DG\rangle =1\,\,\mbox{in}\,\,N_\gamma.
\end{equation*}
The following lemma was proved in \cite{FS:MR833404}.

\begin{Lemma}[{\cite[Lemma~4.1]{FS:MR833404}}]
\label{lem:estwepsilon}
Let $\beta_\eps\to b_0$ uniformly on compact subsets of $N$ as $\eps\to 0$. 
Suppose, for each $\eps>0$, functions  $w^\eps$, $A^\eps$ and numbers $C^\eps$, $a^\eps$  
satisfy
\begin{equation*}
\begin{cases}\ds
-\frac{\eps^2}{2}\Delta w^\eps(x)+\langle \beta_\eps(x),Dw^\eps(x)\rangle =A^\eps(x),\quad x\in N_\gamma,
\\ 
w^\eps(x)=0,\quad x\in N_\gamma\cap\partial\Domain,
\\
|w^\eps(x)|\leq C^\eps,\quad x\in \partial N_\gamma \setminus\partial\Domain,
\end{cases}
\end{equation*}
and $|A^\eps(x)|\leq a^\eps$ for all $x\in N_\gamma$. Then there exists $\eps_0>0$
such that for all $\eps<\eps_0$,
\begin{equation*}
|w^\eps(x)|\leq C^\eps e^{-\frac{2F(x)}{\eps^2}}+4a^\eps G(x),\quad  x\in N_\gamma.
\end{equation*}
\end{Lemma}

\begin{Theorem}
\label{thm:convofvepsilon} As $\eps\to 0$,
\begin{equation*}
v^\eps=v^0+\eps^2 v_1+o(\eps^2)
\end{equation*}
uniformly on compact subsets of $N$. Here, the function $v_1\in C^\infty(N)$ 
is a unique solution of
\begin{equation}\label{d8}
\begin{cases}
\langle b_0(x),Dv_1(x)\rangle =\frac{1}{2}\Delta v^0(x),\quad x\in N,
\\
v_1(x)=0,\quad x\in N\cap\partial\Domain.
\end{cases}
\end{equation}
\end{Theorem}
\bpf
First we recall that since the characteristic curves for (\ref{d8}) are the same
as for (\ref{a3}) in $N$,
$v_1$ is obtained by the method of characteristics in $N$ and is
in $C^\infty(N)$.

Let $N'$ be a compact subset of $N$. It can be covered by a finite number of sets
$N^{\delta_i}_{\gamma_i}=\{x\in N_{\gamma_i}: \dist(x, \partial N_{\gamma_i}
\setminus\partial\Domain)  >\delta_i\}, i=1,...,m$. Let us define
\[
 v^\eps_1=\frac{v^\eps-v^0}{\eps^2}.
\]
We will show that $v^\eps_1\to v_1$ uniformly on every $N^{\delta_i}_{\gamma_i}$.
On $N$ we have
\[
 -\frac{\eps^2}{2}\Delta v^\eps
-\langle b,Dv^\eps\rangle+\frac{1}{2}|Dv^\eps|^2=0
\]
and
\[
 -\frac{\eps^2}{2}\Delta v^0
-\langle b,Dv^0\rangle+\frac{1}{2}|Dv^0|^2
=-\frac{\eps^2}{2}\Delta v^0.
\]
Therefore,
\[
 -\frac{\eps^2}{2}\Delta (v^\eps-v^0)
-\langle b,D(v^\eps-v^0)\rangle
+\frac{1}{2}\langle Dv^\eps+Dv^0,D(v^\eps-v^0)\rangle=
\frac{\eps^2}{2}\Delta v^0,
\]
which gives
\[
 -\frac{\eps^2}{2}\Delta v^\eps_1
+\langle \beta_\eps,Dv^\eps_1\rangle
=\frac{1}{2}\Delta v^0, 
\]
where $\beta_\eps(x)=-b(x)+(Dv^\eps(x)+Dv^0(x))/2$. Lemma~\ref{lem:convofderivatives} implies that
$\beta_\eps\to b_0$ uniformly on compact subsets of $N$.
It follows from (\ref{d8}) that
\[
  -\frac{\eps^2}{2}\Delta v_1
+\langle \beta_\eps,Dv_1\rangle
=\frac{1}{2}\Delta v^0-\frac{\eps^2}{2}\Delta v_1+\langle \beta_\eps-b_0,Dv_1\rangle.
\]
Therefore, the function $w^\eps=v^\eps_1-v_1$ satisfies
\[
 -\frac{\eps^2}{2}\Delta w^\eps
+\langle \beta_\eps,Dw^\eps\rangle
=\frac{\eps^2}{2}\Delta v_1-\langle \beta_\eps-b_0,Dv_1\rangle\quad\mbox{in}\,\,
N_{\gamma_i}
\]
and $w^\eps=0$ on $N_{\gamma_i}\cap\partial\Domain$. Moreover, there is a constant $C$, such that
\[
 |w^\eps|\leq\frac{C}{\eps^2}\quad\mbox{on}\,\,
\partial N_{\gamma_i}\setminus\partial\Domain,\,\,i=1,...,m,
\]
and if $A^\eps(x)=\frac{\eps^2}{2}\Delta v_1(x)-\langle \beta_\eps(x)-b_0(x),Dv_1(x)\rangle$,
then
\[
 |A^\eps(x)|\leq a^\eps\to 0\quad\mbox{on}\,\,\bigcup_{i=1}^m N_{\gamma_i}.
\]
By Lemma \ref{lem:estwepsilon} we thus obtain that for any $i=1,...,m$,
\[
 |w^\eps(x)|\leq \frac{C}{\eps^2} e^{-\frac{2F_i(x)}{\eps^2}}
+4a^\eps G_i(x),
\quad x\in N_{\gamma_i}.
\]
This implies that $w^\eps\to 0$ as $\eps\to 0$ uniformly on 
$N^{\delta_i}_{\gamma_i}$ for any  $i=1,...,m$.
\epf

\subsection{Asymptotics of $Dv^\eps$}

We use the strategy from \cite{F1:MR0304045} where the asymptotics of derivatives
was proved for parabolic problems with zero boundary value on the parabolic
boundary. However we use simpler PDE arguments whenever possible.

\begin{Lemma}
\label{lem:convderv1eps}
$Dv_1^\eps\to Dv_1$ uniformly on compact subsets of $N\cap\partial \Domain$.
\end{Lemma}

\bpf
Let $A$ be a compact subset of $N\cap\partial \Domain$. Since
$v_1^\eps=v_1=0$ on $N\cap\partial \Domain$, we only need to show
$\partial_n v_1^\eps\to \partial_n v_1$ on $A$. Let $x_0=(x_0',0)\in A$. We have
\begin{equation}\label{e1}
\begin{split}
&
\partial_n v_1^\eps(x_0)- \partial_n v_1(x_0)
=
\lim_{x_n\to 0}\frac{v_1^\eps(x_0',x_n)- v_1(x_0',x_n)}{x_n}
\\
&
=\partial_n v^0(x_0)
\lim_{x_n\to 0}\frac{v_1^\eps(x_0',x_n)- v_1(x_0',x_n)}{v^0(x_0',x_n)}.
\end{split}
\end{equation}
Let $0<\sqrt{2} r<{\rm dist}(A,\partial N\cap\Domain)$ and $W_r(x_0)=B_r(x_0')\times[0,r)$, where
$B_r(x_0')$ stands for the Euclidean ball of radius $r$ centered at $x_0'$. We also require
that $r$ is small enough so that there is a constant $c_1>0$ such that 
\begin{equation}\label{eq:e2}
|Dv^0(x)|\geq c_1,\quad x\in W_r(x_0), 
\end{equation}
for every $x_0\in A$.
The functions $w^\eps=\pm(v^\eps_1-v_1)$ satisfy 
\[
 -\frac{\eps^2}{2}\Delta w^\eps(x)
+\langle \beta_\eps(x),Dw^\eps(x)\rangle
=\pm A^\eps(x),\quad x\in W_r(x_0),
\]
where $A^\eps(x)=\frac{\eps^2}{2}\Delta v_1(x)-\langle \beta_\eps(x)-b^0(x),Dv_1(x)\rangle$.
Also, $w^\eps(x)=0$ for $x\in W\cap\partial\Domain$, and there are numbers $a^\eps,c^\eps\to 0$
independent of $x_0\in A$ such that $|A(x)|\leq a^\eps$ for $x\in W_r(x_0)$ and  
$|w^\eps(x)|\leq c^\eps$ for $x\in\partial W_r(x_0)\setminus\partial\Domain$. We also recall that
$\beta_\eps\to b_0$ uniformly on compact subsets of $N$.

Since
\[
\langle b_0(x),Dv^0(x)\rangle=\frac{1}{2}|Dv^0(x)|^2\geq \frac{c_1^2}{2},\quad x\in W_r(x_0), 
\]
it follows that for sufficiently small $\eps$,
\[
\langle \beta_\eps(x),Dv^0(x)\rangle\geq \frac{c_1^2}{4},\quad x\in W_r(x_0). 
\]
 Let $\eta>0$. We set
\[
 \psi(x)=\eta v^0(x)+\frac{c^\eps}{r^2}|x'-x_0'|^2.
\]
Then for $\eps$ small enough (but independent of $x_0$), $\psi\geq c^\eps$ on $\partial W_r(x_0)\setminus\partial\Domain$, and thus $\psi
\geq|w^\eps|$ on $\partial W_r(x_0)$. Moreover, for some constants $c_2,c_3$
\[
 -\frac{\eps^2}{2}\Delta \psi(x)
+\langle \beta_\eps(x),D\psi(x)\rangle\geq -\eps c_2+\frac{\eta c_1^2}{4}
-c^\eps c_3\ge \frac{\eta c_1^2}{8}\ge a^\eps.
\]
if $\eps$ is small enough. Therefore, by comparison we obtain $|w^\eps|
=\max(\pm w^\eps)\leq \psi$ in $W_r(x_0)$. Hence
\[
 \left|\frac{v_1^\eps(x_0',x_n)- v_1(x_0',x_n)}{v^0(x_0',x_n)}\right|
\leq\eta
\]
if $0<x_n<r$, and the claim follows since $\eta$ is arbitrary.
\epf

\bigskip

\bpf[Proof of Theorem~\ref{thm:asymptotics-for-Dv-eps}] Our arguments follow those
of the proof of Theorem~6.4 of~\cite{F1:MR0304045}.

It is sufficient to show that the convergence is uniform on every relatively open subregion of strong regularity $N_1\subset N$ 
such that $\bar {N_1}\subset N$.
Let $N_2$ be a 
relatively open subset of $N$ such that $\bar {N_1}\subset
\bar {N_2}\subset N$.  Let us introduce $b_\eps(x)=-b(x)+Dv^\eps(x)$ so that $\bar b_\eps(x)=-b_\eps(x)$,  see~\eqref{eq:conditioned-drift-eps}. 
We extend $b_0$ and~$b_\eps$ outside $N_2$ to be Lipschitz functions on~$\mathbb{R}^n$
such that
\begin{equation}\label{eq:e3}
 \sup_{\mathbb{R}^n}|b_\eps-b_0|=c^\eps\to 0\quad\mbox{as}\,\,\eps\to 0.
\end{equation}
We will denote by $L$ the Lipschitz constant of $b_0$. For $x\in N_1$, let $X_0$ be the solution of
\begin{equation}\label{eq:e4}
\begin{cases}\ds
\dot X_0(t)=-b_0(X_0(t)),
\\
X_0(0)=x,
\end{cases}
\end{equation}
and $X_\eps$
be the strong solution of the It\^o equation
\begin{equation}\label{eq:e5}
\begin{cases}
dX_\eps(t)=-b_\eps(X_\eps(t))+\eps dW(t),\\
X_\eps(0)=x,
\end{cases}
\end{equation}
where $W$ is a standard $n$-dimensional Wiener process defined on a probability space $(\Omega,\Fc,\Pp)$.
We denote by $\tau^0_x,\tau^\eps_x$ to be respectively the first exit times 
of $X_0$ and $X_\eps$ from $N_2$. We recall that $X_0(t)$ is the characteristic
of (\ref{a3}) and (\ref{d8}) passing through~$x$ with its time parametrization reversed.
Let $T_1$ be such that $\tau_x^0\leq T_1$ for all $x\in N_1$ and let $T=T_1+1$.
Since $\langle b_0(x),e_n\rangle>c_0>0$ for some $c_0$ and all $x\in N_1$
sufficiently close to $\partial\Domain$, there are $c_1,s_0>0$ such that
\begin{equation}\label{eq:e6}
{\rm dist}(X_0(\tau_x^0+s),\partial\Domain)\geq c_1|s|\quad\mbox{for all}\,\,x\in N_1\,\,
\mbox{and}\,\,|s|\leq s_0.
\end{equation}

Let us introduce
$
 A_\eps:=\{\omega\in\Omega:\sqrt\eps\sup_{0\leq t\leq T}|W(t)(\omega)|\leq 1\}$.
Notice that the sets form a monotone family: $A_{\eps_1}\supset A_{\eps_2}$ if $\eps_1\le \eps_2$. Also, $\Pp(\bigcup_{\eps>0}A_\eps)=1$.
A~standard maximal inequality for $W$ implies that for some $C,\gamma>0$,
\begin{equation}\label{eq:e7}
\Pp(A_\eps^c)\leq Ce^{-\frac{\gamma}{\eps}}.
\end{equation}
 It follows from (\ref{eq:e3}), (\ref{eq:e4}), and (\ref{eq:e5})
that
\[
 |X_\eps(t)-X_0(t)|\leq L\int_0^t|X_\eps(s)-X_0(s)|ds+c^\eps t
+\eps |W(t)|.
\]
Therefore, for all $x\in N_1$, $\omega\in A_\eps$,
\[
 |X_\eps(t)-X_0(t)|\leq L\int_0^t|X_\eps(s)-X_0(s)|ds+c^\eps t
+\sqrt\eps,
\]
and, by Gronwall's inequality,
\begin{equation}\label{eq:e8}
 \sup_{0\leq t\leq T}|X_\eps(t)-X_0(t)|\leq k^\eps,
\end{equation}
where $k^\eps= (c^\eps T+\sqrt\eps)e^{TL}\to 0$, as $\eps\to 0$.
Since $X_0(t)\in N_1$ for $0\leq t\leq \tau_x^0$, this implies $X_\eps(t)\in N_2$ for $0\leq t\leq \tau_x^0\wedge \tau_x^\eps$
for $\omega\in A_\eps$
if $\eps<\eps_0$ for some $\eps_0$ independent of $x$ and only depending
on $\dist(N_1,\Domain\setminus N_2)$.
Therefore,  (\ref{eq:e6}) and (\ref{eq:e8}) imply that if $\eps<\eps_0$, $x\in N_1$, and $\omega\in A_\eps$, then
\begin{equation}\label{eq:e9}
 |\tau^\eps_x-\tau_x^0|\leq\frac{k^\eps}{c_1}
\end{equation}
and
$X_\eps(\tau^\eps_x)\in \Gamma$. 

\bigskip

Differentiating equations~(\ref{a2}) and~(\ref{a3})
with respect to $x_i$, $i=1,...,n$, we obtain
\[
 -\frac{\eps^2}{2}\Delta \partial_i v^\eps
+\langle Dv^\eps- b,D\partial_i v^\eps\rangle-\langle \partial_ib,Dv^\eps\rangle=0
\]
and
\[
 -\frac{\eps^2}{2}\Delta \partial_i v^0 
-\langle b,D \partial_i v^0 \rangle-\langle \partial_i b,Dv^0\rangle+\langle Dv^0,D\partial_iv^0\rangle=
 -\frac{\eps^2}{2}\Delta \partial_i v^0.
\]
Subtracting the above equations
and dividing by $\eps^2$ yields 
\begin{multline*}
 -\frac{\eps^2}{2}\Delta \partial_i v^\eps_1
+\langle  Dv^\eps- b,D\partial_i v^\eps_1\rangle
+\frac{1}{\eps^2}\langle  Dv^\eps,D\partial_i v^0\rangle
\\
-\langle \partial_i b,Dv^\eps_1\rangle
-\frac{1}{\eps^2}\langle Dv^0,D\partial_i v^0\rangle
=\frac{1}{2}\Delta \partial_iv^0.
\end{multline*}
Using
\[
\frac{1}{\eps^2}\langle  Dv^\eps,D\partial_i v^0\rangle
-\frac{1}{\eps^2}\langle Dv^0(x),D\partial_i v^0\rangle=\langle D \partial_i v^0,Dv^\eps_1\rangle,
\]
and $\partial_i b_0= -\partial_i b+D\partial_iv^0$, we thus obtain
\[
 -\frac{\eps^2}{2}\Delta \partial_i v^\eps_1
+\langle  b_\eps,D\partial_i v^\eps_1\rangle
+\langle \partial_i b_0,Dv^\eps_1\rangle
=\frac{1}{2}\Delta \partial_i v^0,\quad i=1,\ldots,n.
\]
We can combine these $n$ identities into one:
\begin{equation}\label{e11a}
 -\frac{\eps^2}{2}\Delta (Dv^\eps_1)
+D(Dv^\eps_1)b_\eps
+(Db_0)^*Dv^\eps_1
=\frac{1}{2}\Delta (Dv^0).
\end{equation}
For $\eps\geq 0$, we define the fundamental matrices $Y^\eps$ to be the solutions of
\[
\begin{cases}
 \dot Y^\eps(t)=-Y^\eps(t)(Db_0)^*(X_\eps(t)),
\\
Y^\eps(0)=I.
\end{cases}
\]
Let us denote $\tilde\tau_x^\eps=\tau_x^\eps\wedge T$. Using It\^o's formula and (\ref{e11a}) we obtain
\[
\begin{split}
&
\E[Y^\eps(\tilde\tau_x^\eps)Dv^\eps_1(X_\eps(\tilde\tau_x^\eps))]
\\
&
=Dv^\eps_1(x)+\E\bigg[\int_0^{\tilde\tau_x^\eps}\bigg(-Y^\eps(t)(Db_0)^*(X_\eps(t))Dv^\eps_1(X_\eps(t))
\\
&
\quad\quad
-Y^\eps(t)D(Dv^\eps_1)(X_\eps(t))b_\eps(X_\eps(t))
+\frac{\eps^2}{2}Y^\eps(t)\Delta (Dv^\eps_1)(X_\eps(t))\bigg) dt\bigg]
\\
&
=
Dv^\eps_1(x)-\frac{1}{2}\E\bigg[\int_0^{\tilde\tau_x^\eps}Y^\eps(t)\Delta (Dv^0)(X_\eps(t))dt\bigg],
\end{split}
\]
which yields
\begin{equation}\label{e11}
Dv^\eps_1(x)=\frac{1}{2}\E\left[\int_0^{\tilde\tau_x^\eps}Y^\eps(t)D(\Delta v^0)(X_\eps(t))dt
\right]+
\E[Y^\eps(\tilde\tau_x^\eps)Dv^\eps_1(X_\eps(\tilde\tau_x^\eps))].
\end{equation}
We will show that the right hand side of (\ref{e11}) converges to
\begin{equation}
\label{eq:V}
V(x):=\frac{1}{2}\int_0^{\tau_x^0}Y^0(t)D(\Delta v^0)(X_0(t))dt
+
Y^0(\tau_x^0)Dv_1(X_0(\tau_x^0))
\end{equation}
uniformly on $N_1$. To that end, we need to estimate the difference between the corresponding terms of \eqref{e11} and~\eqref{eq:V}.

We begin with the non-integral terms.
We notice that there is a number $C(T)>0$ such that for every $\omega,\eps$ and $x\in N_1$
\begin{equation}\label{e12}
\sup_{0\leq t\leq \tilde\tau_x^\eps}\|Y^\eps(t)\|\leq C(T).
\end{equation}
Next, there is $\eps_0>0$ such that if $\eps\in[0,\eps_0)$, $\omega\in A_\eps$, and $x\in N_1$,
then (i) $\tilde\tau_x^\eps=\tau_x^\eps$ and (ii) $X_\eps(\tilde\tau^\eps_x)\in \Gamma$. 
Property (i) along with~(\ref{eq:e8}) and standard ODE theory implies that 
there are positive numbers $(k^\eps_1)_{\eps\in(0,\eps_0)}$ such that $k_1^\eps\to 0$ as $\eps\to 0$, and
\begin{equation}\label{e13}
\sup_{0\leq t\leq \tilde\tau_x^\eps}\|Y^\eps(t)-Y^0(t)\|
\leq k_1^\eps,\quad\omega\in A_\eps, x\in N_1.
\end{equation}
Property~(ii) allows us to apply Lemma~\ref{lem:convderv1eps}. So, along with 
(\ref{eq:e8}), (\ref{eq:e9}), (\ref{e12}), and~(\ref{e13}), it implies
that 
there are positive numbers $(k^\eps_2)_{\eps\in(0,\eps_0)}$ such that $k_2^\eps\to 0$ as $\eps\to 0$, and
\begin{equation}\label{e14}
|Y^\eps(\tilde\tau_x^\eps)Dv^\eps_1(X_\eps(\tilde\tau_x^\eps))-Y^0(\tau_x^0)Dv_1(X_0(\tau_x^0))|\leq
k_2^\eps,\quad \omega\in A_\eps, x\in N_1.
\end{equation}
Finally, we observe that there is a constant $C_1$, such that for all $\eps$,
\begin{equation}\label{e15}
|Dv^\eps_1(x)|\leq \frac{C_1}{\eps^2},\quad x\in N_2.
\end{equation}

The difference between the second terms of \eqref{e11} and~\eqref{eq:V} can be 
estimated, due to~\eqref{e12},\eqref{e14},\eqref{e15} as
\begin{align*}
\left|\E[Y^\eps(\tilde\tau_x^\eps)Dv^\eps_1(X_\eps(\tilde\tau_x^\eps))]-
Y^0(\tau_x^0)Dv_1(X_0(\tau_x^0))\right|\le k_2^\eps+C_2(1+\eps^{-2})\Pp(A_\eps^c) 
\end{align*}
for some constant $C_2$, and by \eqref{eq:e7} converges to $0$ as $\eps\to 0$.

To estimate the difference between the integral terms of \eqref{e11} and~\eqref{eq:V}, we write
\begin{align*}
 &\left|\E
\int_0^{\tilde\tau_x^\eps}Y^\eps(t)D(\Delta v^0)(X_\eps(t))dt
-\int_0^{\tau_x^0}Y^0(t)D(\Delta v^0)(X_0(t))dt\right|
\\
\le &\ \E\left[
{\bf 1}_{A_\eps}\int_0^{\tau_x^\eps\wedge\tau_x^0}\left|
Y^\eps(t)D(\Delta v^0)(X_\eps(t))-Y^0(t)D(\Delta v^0)(X_0(t))\right|dt\right] 
\\
&+\E\left[
{\bf 1}_{A_\eps}\left(\int_{\tau_x^\eps\wedge\tau_x^0}^{\tau_x^\eps}\
\left|Y^\eps(t)D(\Delta v^0)(X_\eps(t))\right| dt
+
\int_{\tau_x^\eps\wedge\tau_x^0}^{\tau_x^0} \left|Y^0(t)D(\Delta v^0)(X_0(t))\right|dt\right)\right]
\\
&+\E\left[
{\bf 1}_{A_\eps^c}\left(\int_0^{\tau_x^\eps}\left|
Y^\eps(t)D(\Delta v^0)(X_\eps(t))\right|dt+\int_0^{\tau_x^0} \left|Y^0(t)D(\Delta v^0)(X_0(t))\right|dt\right)\right].
\end{align*}
Each of the terms on the r.h.s.\ uniformly converges to 0 in $N_1$. For the first term, this
follows from~\eqref{eq:e8},\eqref{e12}, and~\eqref{e13}, for the second one ---
from~\eqref{eq:e9} and~\eqref{e12}, and for the third one --- 
from~\eqref{eq:e7} and~\eqref{e12}.

Thus our claim of uniform convergence of  $Dv^\eps_1$ to $V$ follows. 
It remains to notice, differentiating~(\ref{d8}), that $Dv_1$ satisfies in $N$ the system
\[
\langle b_0,D \partial_i v_1 \rangle+\langle \partial_i b_0,Dv_1(x)\rangle
=\frac{1}{2}\Delta \partial_i v^0, \quad i=1,...,n,
\]
for which the method of characteristics implies that $Dv_1(x)=V(x)$ in $N_1$.
\epf

\section{Doob's $h$-transform for conditioned diffusions}\label{sec:h-transform}

Here we provide a general and rigorous introduction to Doob's $h$-transform computing
the conditional distribution for diffusions conditioned on exit events. 
The material of this section is not highly original: 
connections of $h$-transform to conditioning, to the potential theory of elliptic PDEs, and 
to Martin boundaries are  well-known, see, e.g.,~\cite{Doob:MR731258,Pinsky:MR1326606}. However, 
no rigorous and complete exposition of the
main result of this section, Theorem~\ref{thm:h-transform-for-diffusions}, 
is known to us, and we decided 
to include this section, hoping that it will serve as a useful reference point for future research.
Our exposition is based on \cite{Bloemendal} and~\cite[Chapter~13]{Wentzell:MR781738}. We use minimal information from the PDE theory. 
Other useful sources on Markov processes and diffusions are~\cite{Stroock--Varadhan:MR532498}, \cite{Ethier-Kurtz:MR838085},
\cite{Ikeda-Watanabe:MR1011252}, \cite{Karatzas-Shreve:MR1121940}.

Let us first introduce an abstract generalization of a diffusion process in a domain with absorption at the boundary of the domain.
We will always work with homogeneous Markov processes, i.e., processes with transition mechanisms that do not depend on initial time.

Let $\Domain$ be a domain in $\R^n$. Let us equip the space $C=C([0,\infty),\bar \Domain)$ of continuous paths 
$(X(t))_{t\ge 0}$ with $\Bc=\Bc(C)$, the Borel $\sigma$-algebra with respect to locally uniform topology.

Suppose that $\Pp^t(x,dy)$ is a Markov transition kernel on $\bar \Domain$. It means that (i) for all $t\ge 0$ 
and all $x\in\bar \Domain$, $\Pp^t(x,\cdot)$ is a Borel probability measure on~$\bar \Domain$; (ii) for any $t\ge 0$  and any Borel set $A$,
$\Pp^t(\cdot,A)$ is a Borel measurable function; (iii)  $\Pp^0(x,dy)=\delta_x(dy)$; (iv) the Chapman--Kolmogorov equations hold, i.e.,
for any $s,t\ge0$ and any Borel set $A$, 
\[
 \Pp^{t+s}(x,A)=\int_{\bar \Domain}\Pp^t(x,dy)\Pp^s(y,A).
\]
Let us make the following additional assumptions:
\begin{enumerate}
\item \label{prop:Markov-family} There is a family of measures $(\Pp_x)_{x\in\bar\Domain}$ on paths $(C,\Bc)$ 
such that for each $x\in\bar \Domain$,
$\Pp_x\{X(0)=x\}=1$ and under $\Pp_x$ 
the process $X$ is a Markov process with a homogeneous transition probability $\Pp^t(x,dy)$. 
\item \label{prop:preserves-boundary}
For all $x\in\partial \Domain$ 
and all $t\ge 0$, $\Pp^t(x,dy)=\delta_x(dy)$.
\end{enumerate}

Processes associated with such transition kernels or semigroups will be called
continuous Markov processes on~$\Domain$ with absorption at~$\partial\Domain$. 
We denote by~$\E_x$ the expectation with respect to $\Pp_x$. The semigroup $(\Pp^t)$
is defined by
\[
 \Pp^t f(x)=\int_{\bar \Domain} \Pp^t(x,dy)f(y),\quad t\ge 0,\ x\in\bar \Domain,\ f\in\Bb(\bar\Domain), 
\]
where $f\in\Bb(\bar\Domain)$ is the space of bounded measurable 
functions on $\bar\Domain$.

For any $X\in C$,  we denote by $\tau(X)$ the first exit on the boundary: $\tau=\inf\{t\ge 0:\ X(t)\in \partial \Domain\}$.

Let $\Gamma$ be a measurable subset of $\partial \Domain$. We introduce a trajectory set
\[
C_\Gamma=\{X\in C:\ \tau(X)<\infty,\ X(\tau(X))\in\Gamma\}
\]
and a measurable bounded function
\begin{equation}
\label{eq:definition_of_h}
h(x)=h_\Gamma(x)=\Pp_x(C_\Gamma)=\lim_{n\to\infty}P^n(x,\Gamma),\quad x\in\bar\Domain.
\end{equation}
Let us assume that 
\begin{equation}
\label{eq:h-positive}
h(x)>0, \quad x\in\Domain.
\end{equation}
Our goal is then to describe the conditional measures  $\Pp_{\Gamma,x}$ defined  by
\[
 \Pp_{\Gamma,x}(A)= \Pp_x(A| C_\Gamma),\quad x\in\Domain\cup \Gamma, A\in\Bc.
\]
We will denote expectation with respect to $\Pp_{\Gamma,x}$ by $\E_{\Gamma,x}$.

Denoting by $\Fc_t$ the natural filtration of the process $X$,
we obtain for any $\Fc_t$-measurable random variable $\xi$:
\begin{equation}
\label{eq:harmonic-flow-of-measures}
\E_{\Gamma,x}\xi=\E_x\left[\xi \frac{\ONE_{C_\Gamma}}{h(x)}\right]=
\E_x\left[\xi\E_x\left[\frac{\ONE_{C_\Gamma}}{h(x)}\Bigr|\Fc_t\right]\right]
=\E_x\left[\xi\frac{h(X_t)}{h(x)}\right].
\end{equation}

\begin{Lemma} If $x\in\Domain\cup\Gamma$, then $\Pp_{\Gamma,x}$ defines a continuous
Markov process on $\Domain\cup\Gamma$ with transition probability
\begin{equation}
\label{eq:P-Gamma}
\Pp_{\Gamma}^t(x,dy)=\frac{h(y)}{h(x)}\Pp^t(x,dy).
\end{equation}
\end{Lemma}
\bpf
The continuity is inherited from the original process, so 
it is sufficient to show that for any bounded measurable function $f$ on $\Domain\cup\Gamma$,
and any $s,t\ge 0$, $x\in \Domain\cup\Gamma$,
\begin{equation*}
\E_{\Gamma,x}[f(X(s+t))|\Fc_s]=\E_x\left[f(X(s+t))\frac{h(X(s+t))}{h(X(s))}\Bigr|\ X(s)\right].
\end{equation*}
The right-hand side is $\Fc_s$-measurable, so we need to 
check that integrals of both sides over any event $A\in\Fc_s$ coincide. By~\eqref{eq:harmonic-flow-of-measures},
the integral identity to check is
\begin{multline*}
\E_x\left[\E_x\left[f(X(s+t))\frac{h(X(s+t))}{h(X(s))}\Bigr|\ X(s)\right]\ONE_{A}\frac{h(X(s))}{h(x)}\right]
\\=\E_x\left[f(X(s+t))\ONE_{A}\frac{h(X(s+t))}{h(x)}\right].
\end{multline*}
To prove this identity we cancel the two instances of $h(X(s))$ on the left-hand side and, using Markov property of $X$ under $\Pp_x$, replace 
conditioning with respect to $X(s)$ by conditioning
with respect to $\Fc_s$.\epf

If $x\in\partial\Domain\setminus\Gamma$, then the above construction does not make sense, and we simply set
$\Pp^t_\Gamma(x,dy)=\delta_x(dy)$ for all $t\ge 0$. 
Combining this with~\eqref{eq:P-Gamma} we see that thus defined process is also a continuous Markov process in $\Domain$ absorbed at $\partial\Domain$,
and the  action of the semigroup associated with transition kernels~$\Pp^t_\Gamma$
can be written as
\begin{align}
\label{eq:P-Gamma-semigroup}
\Pp^t_\Gamma f(x)&=\frac{\E_x f(X(t))h(X(t))}{h(x)}\ONE_{x\in \Domain}+f(x)\ONE_{x\in\partial\Domain}\\
\notag
&=\frac{\Pp^t (hf)(x)}{h(x)} \ONE_{x\in \Domain} +f(x)\ONE_{x\in\partial\Domain},\quad t\ge 0,\ x\in\bar \Domain,\   f\in \Bb(\bar\Domain).
\end{align}

\medskip

Let us say that a semigroup $(\Pp^t)$  defines 
a diffusion process in $\Domain$ absorbed at $\partial\Domain$
if in addition to properties~\ref{prop:Markov-family} and~\ref{prop:preserves-boundary} the following holds:
for each point $x\in\bar\Domain$ there is 
a  positive semi-definite symmetric $n\times n$ matrix
$a(x)$ and an $n$-dimensional vector $b(x)$ such  that $a$ and $b$ are Borel functions on $\bar\Domain$, bounded
on every compact subset of $\Domain$, and for every function $f\in C^2_0(\Domain)$
(i.e.,  $f\in C^2(\Domain)\cap C^0(\bar\Domain)$ and $\supp f \subset \Domain$),
the generator
\[
Af=\lim_{t\to 0}\frac{\Pp^tf -f}{t}
\]
is well defined in the space $\Bb(\bar\Domain)$ equipped with sup-norm (i.e., the convergence in the right-hand side is uniform) and 
\begin{equation}
\label{eq:generator-equals-second-order-op}
Af(x)=Lf(x),\quad  x\in\Domain,
\end{equation}
where we denote
\begin{equation}
\label{eq:second-order-operator}
Lf(x)= \sum_{i=1}^n b^i(x)\partial_i f(x)+
\frac{1}{2}\sum_{i,j=1}^n a^{ij}(x)\partial_{ij}f(x),
\end{equation}
whenever the derivatives involved are well-defined.

Notice that we require~\eqref{eq:generator-equals-second-order-op} to hold for $x\in\Domain$, although it is often
convenient to have coefficients $a$ and $b$ defined on $\bar\Domain$.

The identification of diffusion processes with solutions of stochastic equations and martingale problems is well-known,
see, e.g.,~\cite{Stroock--Varadhan:MR532498}, \cite{Ethier-Kurtz:MR838085}, \cite{Karatzas-Shreve:MR1121940}. 
The following theorem claims that diffusion processes with absorption can also be viewed as solutions of It\^o stochastic equations
stopped upon reaching the boundary.

\begin{Theorem}\label{thm:sde-representation} Suppose the semigroup $(\Pp^t)$ defines a diffusion process on a domain $\Domain$ with absorption at $\partial\Domain$ and
coefficients $a,b$. Let $\sigma$ be a Borel measurable $n\times n$ matrix-valued function on $\Domain$ such that $\sigma\sigma^*\equiv a$
and such that~	$\sigma$ is bounded on any compact subset of $\Domain$.
Then for any $x\in\Domain$ there is an extension  $(\tilde C,\tilde \Bc,\tilde \Pp_x)$  of the original probability space $(C,\Bc,\Pp_x)$,
a filtration $(\tilde\Fc_t)_{t\ge 0}$ on $\tilde\Bc$
satisfying the usual conditions, and an $n$-dimensional 
Wiener process~$W$ w.r.t. $(\tilde \Fc_t)$, such that the coordinate process $X$ is $(\tilde \Fc_t)$-adapted and,
with probability one, $X(0)=x$ and 
\[
 dX(t)=b(X(t))\ONE_{\{X(t)\in\Domain\}}dt+\sigma(X(t))\ONE_{\{X(t)\in\Domain\}}dW(t),\quad t\ge 0.
\] 
If $a(x)$ is nondegenerate for all $x\in\Domain$, then $(\tilde C,\tilde \Bc,\tilde \Pp_x)$  may be taken to coincide with $(C,\Bc,\Pp_x)$.
\end{Theorem}

\bpf[Sketch of proof] We only indicate the changes one needs in
adapting the proof of the same statement for the usual diffusion processes, see Proposition~5.4.6 
of~\cite{Karatzas-Shreve:MR1121940}. Let us consider 
a sequence $(U_m)_{m\in\N}$ of open sets such that 
 $K_m=\bar U_m$ is a compact subset of $\Domain$ for every $m\in\N$, and
$\bigcup_{m\in\N} U_m=\Domain$.
For every $m\in\N$, we can find a nonnegative bounded function $g_m\in C^2_0(\Domain)$ such that $g_m\equiv 1$ on $K_m$.
Then functions $f_{m,i}(x)=x_ig_m(x)$ and $f_{m,i,j}(x)=x_ix_jg_m(x)$ belong to the domain of $A$, and~\eqref{eq:generator-equals-second-order-op}
holds for $f=f_{m,i}$ and  $f=f_{m,i,j}$.
Therefore, if~$X$ is the coordinate process under~$\Pp_x$, then, for $f=f_{m,i}$ and  $f=f_{m,i,j}$, the process
$f(X(t))-\int_0^t Lf(X(s))ds$ is a bounded martingale under~$\Pp_x$ (see, e.g.,~\cite[Proposition~4.1.7]{Ethier-Kurtz:MR838085}). 
As in the proof of Proposition~5.4.6 in~\cite{Karatzas-Shreve:MR1121940}, we can use this to derive that
if $\tau_m=\inf\{t\ge 0: X(t)\notin U_m\}$, then
\[
 M_{m}^i(t)=X^i(t\wedge \tau_m)-x_{0}^i-\int_0^{t\wedge \tau_m}b^i(X(s))ds
\]
is a martingale with
\[
 \langle M_{m}^i,M_{m}^j\rangle=\int_0^{t\wedge \tau_m}a^{ij}(X(s))ds.
\]
We can then follow the proof of Proposition~5.4.6 in~\cite{Karatzas-Shreve:MR1121940} and use 
the multi-dimensional version of Doob's representation for continuous martingales (see~\cite[Theorem~3.4.2 and Remark~3.4.3]{Karatzas-Shreve:MR1121940})
to represent $M_m=(M_{m}^i)_{i=1}^n$ as
\[
 M_{m}(t)=\int_0^{t\wedge \tau_m}\sigma(X(s))dW_m(s)
\]
for a Wiener process $W_m$ on a filtered extension $(\tilde C_m,\tilde \Bc_m,(\tilde \Fc_{m,t})_{t\ge 0},\tilde \Pp_{x,m})$ of the original probability space. 
In fact,  one can choose the extended probability space and the Wiener process $W=W_m$ to be independent of $m$, thus obtaining
\[
 X(t\wedge \tau)=x_0+\int_0^{t\wedge \tau}b(X(s))ds+ \int_0^{t\wedge \tau}\sigma(X(s))dW(s),\quad s\ge 0,
\]
where $\tau=\lim_{m\to\infty}\tau_m\in(0,\infty]$ and we used the continuity of trajectories of $X$.
\epf

\bigskip

Theorem~\ref{thm:sde-representation} uses the existence of a Markov process with given coefficients of drift and diffusion as an assumption.
In general, verifying this assumption may be a nontrivial issue. We will not give rigorous general results, but 
let us mention that such results will be parallel to existence
results for the usual diffusion processes and sketch several approaches without implementing them rigorously. 

The first approach is to construct the semigroup directly, using general existence results and the properties of the generator
of diffusion with absorption, see, e.g., Theorem~8.1.4 in~\cite{Ethier-Kurtz:MR838085} for the case of a bounded smooth domain
and H\"older coefficients.

Another approach is to construct the diffusion process $(X(t))_{t\ge 0}$ in $\R^n$, introduce
a stopping time $\tau$ as the hitting time for $\Domain^c$ or $\partial\Domain$, and prove that the stopped
process $(X(t\wedge\tau))_{t\ge 0}$ is a Markov process with all the required properties.
Although this procedure corresponds precisely to the intuition on diffusions stopped upon reaching a closed set, 
the limitation of this method is that it requires smooth continuation of~$a$ and~$b$ through~$\partial O$, 
which may not be possible for irregular boundaries if the coefficients are given only inside $\Domain$.

One more approach is to mimic the proof of Theorem~\ref{thm:sde-representation} and use existence results for SDEs by sequentially constructing solutions
on expanding compact domains. We omit the details of this procedure and notice only that without
imposing some restrictions on the coefficients, such as global Lipschitzness of $a$ and $b$ or existence of an appropriately understood Lyapunov function,
one cannot exclude finite time explosion of solutions.

\bigskip

To derive  the conditioned diffusion generator on $C^2_0(\Domain)$, we will need to check that 
formula~\eqref{eq:generator-equals-second-order-op} holds true for $f=h$. 
The Markov property implies that $h$ is harmonic for $(\Pp^t)$, i.e.,
\[
\Pp^t h(x)= h(x),\quad x\in\bar\Domain,\ t>0,  
\]
Therefore, $Ah$ is well defined and identically equal to $0$ on $\bar\Domain$. However, $h\notin C^2_0(\Domain)$,
and {\it a~priori} it is not even clear if $h\in C^2(\Domain)$. To claim the latter one needs to make certain assumptions.
\begin{Lemma}\label{lem:h-is-C2} Let the coefficients $a,b$ of a diffusion process 
 in $\Domain$ stopped at~$\partial\Domain$ satisfy
$a,b\in C^1(\bar\Domain)$ and $\det a(x)\ne 0$ for all $x\in\bar\Domain$. Then $h\in C^2(\Domain)$.
\end{Lemma}
We postpone the proof of this lemma to the end of this section.

\begin{Lemma}\label{lem:Lh} Let conditions of Lemma~\ref{lem:h-is-C2} hold. Then
\begin{equation}
\label{eq:h-harmonic-generator}
Ah(x)=Lh(x)=0,\quad x\in\Domain.
\end{equation}
\end{Lemma}
\bpf We cannot apply~\eqref{eq:generator-equals-second-order-op} directly since $h\notin C^2_0(\Domain)$.
Let us take $r>0$ such that the ball $B_r(x)$ is contained in $\Domain$. Then $a$ and $b$ are uniformly bounded in $B_r(x)$.
Using Theorem~\ref{thm:sde-representation} to represent the process 
as a solution of an SDE and noticing that the behavior of the process  until the first exit from $B_r(x)$ is entirely determined by
the behavior of the coefficients in $B_r(x)$, one can derive from standard maximal inequalities for martingales that
\[
\Pp_x\left\{\sup_{s\in[0,t]}|X(s)-x|\ge r\right\}=o(t).
\]
In particular, $\Pp^t(x,B_r(x)^c)=o(t)$. Let us now find $f\in C^2_0(\Domain)$ such that
$0\le f(y)\le 1$ for all $y\in\bar\Domain$ and $f(y)=h(y)$ for all $y\in B_r(x)$. Then 
\begin{multline*}
\left|\frac{\Pp^th(x)-h(x)}{t}-\frac{\Pp^t f(x)-f(x)}{t}\right|
= \left|\frac{\int_{\bar\Domain} (h(y)-f(y))\Pp^t(x,dy)}{t}\right|
\\ \le \frac{2\Pp^t(x,B_r(x)^c)}{t}\to 0,\quad t\to 0.
\end{multline*}
Therefore, $0=Ah(x)=Af(x)$ and since the partial derivatives of $h$ and $f$ coincide at $x$, 
formula~\eqref{eq:h-harmonic-generator} is implied by~\eqref{eq:generator-equals-second-order-op} \epf

\begin{Lemma}\label{lem:h-positive} Let conditions of Lemma~\ref{lem:h-is-C2} hold and assume that there is
a point $x_0\in\Domain$ such that $h(x_0)>0$. Then~\eqref{eq:h-positive} holds.
\end{Lemma}
\bpf The lemma is a direct consequence of Lemma~\ref{lem:Lh} and the strong maximum principle, see \cite[Theorem 3, page 349]{Evans:MR2597943}
(or Harnack inequality, see \cite[Corollary 9.25]{Gilbarg-Trudinger:MR1814364}).
One can also give a more probabilistic argument:
Lemma~\ref{lem:h-is-C2} implies that $h$ is positive in some neighborhood~$U$ of~$x_0$. 
For any starting point $x\in\Domain$ there is a continuous path $\gamma:[0,1]\to\Domain$ connecting $x$ to $x_0$,
which implies that $\Pp^1(x, U)>0$. Now $h(x)>0$ follows from the Markov property.
\epf

\begin{Theorem}\label{thm:h-transform-for-diffusions}
 Let $\Pp^t$ define a diffusion process $X$ in a domain $\Domain$, 
with coefficients $a,b\in C^1(\bar\Domain)$ such that $\det a(x)\ne 0$ for all $x\in\bar\Domain$.
Let a measurable set $\Gamma\subset\partial\Domain$ be such that $h(x_0)>0$ for some $x_0\in\Domain$,  where $h$ is defined by~\eqref{eq:definition_of_h}.  Then
the process $X$ conditioned on exit from $\Domain$ through $\Gamma$ is a diffusion process in $\Domain$
with coefficients $a_\Gamma, b_\Gamma$, where $a_\Gamma\equiv a$ and 
\begin{equation*}
b_\Gamma(x)=b(x)+a(x)\frac{Dh(x)}{h(x)},\quad x\in \Domain.
\end{equation*}
\end{Theorem}
\begin{Remark}\rm Under the conditions on $a$ and $b$ imposed by  Theorem~\ref{thm:h-transform-for-diffusions},
the condition $h(x_0)>0$ is effectively a restriction on the ``size'' of $\Gamma$. The domain $\Domain$ itself is not required
to be bounded, and $\partial\Domain$ can be arbitrarily irregular. One natural situation where this condition 
holds true is where $\Gamma$ contains a smooth hypersurface $\Gamma'$ such that $\Gamma'\cup\Domain$ is path-connected.
Also, one can often make sense of $h$-transform for semigroups and their generators even for conditioning on events of zero probability.
\end{Remark}

\bpf[Proof of Theorem~\ref{thm:h-transform-for-diffusions}] 
Lemma~\ref{lem:h-is-C2} and condition~\eqref{eq:h-positive} imply that if $f\in C^2_0(\Domain)$, then $hf\in C^2_0(\Domain)$, and
$h(x)$ is bounded away from zero on the support of $f$.
So,
the generator $A_\Gamma$ of the semigroup $(\Pp^t_\Gamma)$ 
is, by~\eqref{eq:P-Gamma-semigroup}, well-defined on $C^2_0(\Domain)$
and given by
\begin{align*}
A_\Gamma f &=\lim_{t\to 0}\frac{h^{-1}\Pp^t(hf) -f}{t}=
\lim_{t\to 0}\frac{\Pp^t(hf) -hf}{th}=\frac{A(fh)}{h},\quad f\in C^2_0(\Domain),
\end{align*}
and a straightforward computation using Lemma~\ref{lem:Lh} produces
\begin{align*}
A_\Gamma f(x)&=\sum_{i=1}^n \left(b^i(x)+\frac{1}{h(x)}\sum_{j=1}^na^{ij}(x)\partial_j h(x)\right)\partial_i f(x)+
\frac{1}{2}\sum_{i,j=1}^n a^{ij}(x)\partial_{ij}f(x)
\\
&= Lf(x)+\frac{1}{h(x)}\sum_{i,j=1}^n a^{ij}(x)\partial_j h(x)\partial_i f(x),\quad x\in\Domain,
\end{align*}
which completes the proof.
\epf

\bpf[Proof of Lemma~\ref{lem:h-is-C2}]
Let us first summarize the necessary information from the theory of Green functions 
of elliptic PDE's in bounded smooth domains (see, e.g., \cite[Sections~10, 16, 21, 36]{Miranda:MR0087853}).
 
 \begin{Theorem}\label{thm:Green}
Let $B$ and $B'$ be two balls in $\R^n$ and $\bar B\subset B'$. Let
$a,b\in C^1(B')$ for some ball $B\subset B'$ and let $a(x)$ be nondegenerate for all $x\in B'$. Then
there is a function  
$K_B\in C(B\times\partial B)$ such that
for any $\phi\in C(\partial B)$ there is a unique solution
solution $v\in C(\bar B)\cap C^2(B)$ of
\begin{equation*} 
\begin{cases}
Lv(x)=0,& x\in B,
\\
v(x)=\phi(x),& x\in\partial B,
\end{cases}
\end{equation*}
(where $L$ is defined in~\eqref{eq:second-order-operator}) can be represented as 
\begin{equation}
\label{eq:potential-solution}
v(x)=\int_{\partial B}K_{B}(x,y)\phi(y) \mu_B(dy),\quad x\in B
\end{equation}
where $\mu_B(dy)$ denotes the surface area on the sphere $\partial B$.
If $\phi\in C^3(\partial B)$, then $v\in C^2(\bar B)$ and it can be extended to a function $v\in C^2_0(\Domain)$.
 \end{Theorem}

Let us now recall the connection with diffusions.

Consider transition probabilities $(\Pp^t)$ and the associated Markov family $(\Pp^t_x)$  
of diffusion processes with generator $A$ satisfying the conditions of Theorem~\ref{thm:h-transform-for-diffusions}. 
Let $B$ be a ball such that  $\bar B\subset\Domain$ and
let $x\in B$.
Under~$\Pp_x$, we define $\tau_B =\inf\{t\ge 0:\ X(t)\in\partial B\}$. Under the assumptions of Theorem~\ref{thm:h-transform-for-diffusions}, 
$\Pp_x\{\tau_B<\infty\}=1$, see, e.g.,~\cite[Proposition 8.2]{Bass:MR1483890}.

For any $\phi\in C^3(\partial B)$
and the associated function $v$ given by Theorem~\ref{thm:Green} and extended to a function in $C^2_0(\Domain)$ belongs to the domain of
the generator~$A$, and so
$v(X(t))-\int_0^t Av(X(s))ds$ is a bounded martingale under~$\Pp_x$ (see, e.g.,~\cite[Proposition~4.1.7]{Ethier-Kurtz:MR838085}). Since $v\in C^2_0(\Domain)$, 
we have $Av(X(t))\ONE_{t\le \tau_B}=Lv(X(t))\ONE_{t\le \tau_B}=0$, and
Doob's optional sampling theorem~\cite[Theorem 2.2.13]{Ethier-Kurtz:MR838085} implies
\[
v(x)=\E_x\phi(X(\tau_B)),\quad x\in B.
\]
Comparing this to~\eqref{eq:potential-solution}, we conclude that for any 
ball $B\subset \Domain$ and any starting point $x\in\Domain$, 
$K_B(x,\cdot)$ is the density of the distribution of $X(\tau_B)$ with respect to the surface area 
on $\partial B$, so that for any bounded measurable function $f:\partial B\to\R$,
\begin{equation}
\label{eq:density-representation}
\E_x f(X(\tau_B))=\int_{\partial B} f(y) K_B(x,y)\mu_B(dy).
\end{equation}
Under the conditions of Theorem~\ref{thm:h-transform-for-diffusions} the Feller property holds true 
(see, e.g.,~\cite[Theorem 8.1.4]{Ethier-Kurtz:MR838085})
and hence
due to the continuity of trajectories, so does strong Markov property (see, e.g.,~\cite[Theorem 4.2.7]{Ethier-Kurtz:MR838085}). 
The latter implies 
$h(x)=\E_x h(X(\tau_B))$ for any ball $B\subset\Domain$ and any $x\in B$. This, along with~\eqref{eq:density-representation},
implies
\begin{equation*}
h(x)=\int_{\partial B} h(y) K_{B}(x,y)\mu_{B}(dy),\quad x\in B.
\end{equation*}
In the right-hand side, $h$ is bounded and, for any open set $U$
compactly contained in $B$, the function $K_B$ is uniformly continuous on $U\times\partial B$.
Therefore, $h\in C(B)$. Let us take another ball $B'$ such that $\partial B'\subset B$. Then
$h\in C(\partial B')$ and, for $\tau'=\inf\{t\ge0:\ X(t)\in\partial B'\}$
\begin{equation*}
h(x)=\E_x h(X(\tau'))=\int_{\partial B'} h(y) K_{B'}(x,y)\mu_{B'}(dy),\quad x\in B'.
\end{equation*}
Theorem~\ref{thm:Green} implies $h\in C^2(B')$. Since one can choose balls $B$ and $B'$ to contain
any given point in $\Domain$, we conclude that $h\in C^2(\Domain)$. \epf

\bpf[Alternative proof of Lemma~\ref{lem:h-is-C2}]
The proof of Lemma~\ref{lem:h-is-C2} was based on existence of Green's function. We want to present an alternative approach
based on strong Feller property of diffusions (which heuristically means that transition probabilities have nice densities).
Let us recall, see, e.g.,~\cite[Theorem 7.2.4]{Stroock--Varadhan:MR532498}, that if
coefficients $a$ and $b$ are bounded on $\R^n$, $a\in C(\R^n)$, and $\det a>c_0$ for some $c_0>0$,
then the corresponding diffusion process is strong Feller, i.e., for any bounded measurable function $f:\R^n\to \R$ and any time $t>0$,
the function $\Pp^tf$ defined by $\Pp^tf(x)=\int_{\R^n}\Pp^t(x,dy)f(y)$ is continuous. (Another proof of
strong Feller property when the coefficients are Lipschitz continuous can be deduced from the results of 
~\cite[Chapter~5, Sections~4 and~5]{Friedman:MR0494490}.)
So, for any open balls $B_1$ and $B_2$ such that $\bar B_1\subset \bar B_2$ and $\bar B_2\subset\Domain$, we can find 
bounded coefficients $\tilde a,\tilde b\in C^1(\R^n)$  such that 
\begin{equation}
\label{eq:localized-coefs}
 a(x)=\tilde a(x),\quad b(x)=\tilde b(x),\quad x\in\bar B_2,
\end{equation}
and the diffusion associated with $\tilde a,\tilde b$ is strong Feller on $\R^n$.
Let us also extend the function $h$ to a bounded measurable function defined on $\R^n$. 
For any $x\in\bar B_1$, we can use Theorem~\ref{thm:sde-representation} to realize the diffusion 
corresponding to coefficients $a,b$ with absorption at $\partial\Domain$ as a solution to an SDE
driven by a Wiener process.
We can now construct a strong solution of the SDE driven by the same Wiener process, with 
coefficients $\tilde a,\tilde b$ and starting point~$x$
on the same probability space.
Let us keep $\Pp_x$ and $\E_x$
as the notation for the respective probability and expectation on this probability space and denote the diffusion 
processes by $X(t)$ and $\tilde X(t)$. These processes coincide at least up to a random time $\nu=\inf\{t\ge 0: X(t)\in \partial B_2\}$.

We know that the measurable function $h$ satisfies
\[
 h(x)=\E_x h(X(\tilde\tau))
 \]
for every stopping time $\tilde\tau<\tau$.
We need to show that $h$
is continuous in $\Domain$. Denote $\nu_t:=t\wedge\nu$.
Due to~\eqref{eq:localized-coefs}, for any $t\ge 0$,
\begin{equation}
\label{eq:continuous-approx-to-h}
 h(x)=\E_x h(X(\nu_t)) =\alpha_t(x)+\beta_t(x),\quad x\in B_1,
\end{equation}
where $\alpha_t(x)=\E_x h(\tilde X(t))$ and
$\beta_t(x)=\E_x h(X(\nu_t))\ONE_{\{\nu<t\}}-\E_x h(\tilde X(t))\ONE_{\{\nu<t\}}$.
The strong Feller property
for $\tilde X$ implies that $\alpha_t(\cdot)$ is continuous on $B_1$. 
For the second term we have
\[
 |\beta_t(x)|\le 2\Pp_x\{\nu<t\},
\]
and the standard maximal inequalities imply that, as $t\to 0$, $\beta_t(\cdot)$ converges to 0 uniformly in~$B_1$.
Therefore, due to~\eqref{eq:continuous-approx-to-h}, $h$ is continuous on $B_1$ being a uniform limit of continuous functions.
Since the choice of $B_1$ was arbitrary, $h$ is continuous on~$\Domain$.

Once we know that $h$ is continuous in $\Domain$, for every open ball $B$ such that $\bar B\subset\Domain$, the problem
\begin{equation*} 
\begin{cases}
Lv(x)=0,& x\in B,
\\
v(x)=h(x),& x\in\partial B,
\end{cases}
\end{equation*}
has a unique solution $v\in C^2(B)\cap C(\bar B)$ (see e.g. \cite[Theorem 6.13]{Gilbarg-Trudinger:MR1814364}). For any $x\in B$,
we can use the It\^o formula along with the martingale property to see
 \[
 h(x)=\E_x h(X(\tau_B))=\E_x v(X(\tau_B))=v(x)+\E_x\left[\int_0^{\tau_B}Lv(X(t))dt\right]=v(x),
 \] 
 where  $\tau_B =\inf\{t\ge 0:\ X(t)\in\partial B\}$.
 So, $h$ coincides with $v$ in $B$. Therefore, $h\in C^2(B)$. Since the choice of $B$ is arbitrary, the lemma follows.
 \epf

\bibliographystyle{alpha}
\bibliography{exitclt}
\end{document}